\documentclass{article}
\usepackage[russian,english]{babel}
\usepackage[T2A]{fontenc}
\usepackage[cp1251]{inputenc}
\usepackage{amssymb, amsthm, amsmath}
\usepackage{graphicx}
\usepackage{float}
\usepackage[usenames]{color}
\usepackage{colortbl}
\usepackage{subfigure}

\newtheorem{theorem}{Theorem}

\newtheorem{lemma}{Lemma}

\newtheorem{statement}{Statement}

\begin{document}
	\sloppy
	\date{}
	\title{Stable components for gradient-like diffeomorphisms of torus inducing matrix $\begin{pmatrix} -1 & -1\cr 1& 0\end{pmatrix}$} 
	\author{D. Baranov, O. Pochinka}
	\maketitle
	
	\begin{abstract} An isotopy between two  diffeomorphisms $f_0,f_1:M\to M$ means the existence of an arc $\{f_t: M \to M, t \in [0,1]\}$ connecting them in the space of diffeomorphisms. Among such arcs there are so-called stable arcs, which do not qualitatively change under small perturbations. In the present paper we consider a set of gradient-like diffeomorphisms $f$ of 2-torus $\mathbb T^2$ whose induced isomorphism $f_*:\pi_1(\mathbb T^2)\to\pi_1(\mathbb T^2)$ given by a matrix $f_\bigstar=\begin{pmatrix} -1 & -1\cr 1& 0\end{pmatrix}$. We prove that the set of such diffeomorphisms is decomposed into four stable components. Moreover, we establish that two diffeomorphisms under consideration are stably connected  if and only if they have the same number of fixed sinks. 
	\end{abstract}


	
	\section{Introduction and formulation of results}
	Let $M^n$  be a closed smooth connected $n$-manifold. Diffeomorphisms $\varphi_0,\varphi_1:M^n\to M^n$ are called {\it diffeotopic}, if there exists a family of diffeomorphisms $\varphi_t:M^n\to M^n$ (an {\it arc}) that smoothly depends on the parameter $t\in[0,1]$ and is called. 	According to \cite{NPT}, an arc $\varphi_t:M^n\to M^n$ is called {\it stable} if it is an interior point of an equivalence class with respect to the following relation: arcs $\varphi_t$, $\varphi'_t$ are called {\it conjugate} if there exist homeomorphisms $ H:~[0,1]~\to~[0,1 ], \, H_t:~M^n \to M^n$ such that $H_t \varphi_t = \varphi' _ { t} H_t$, and $ H_t $ depends continuously on $t$.
	
\begin{figure}[h!]
		\centering\includegraphics[scale=0.3]{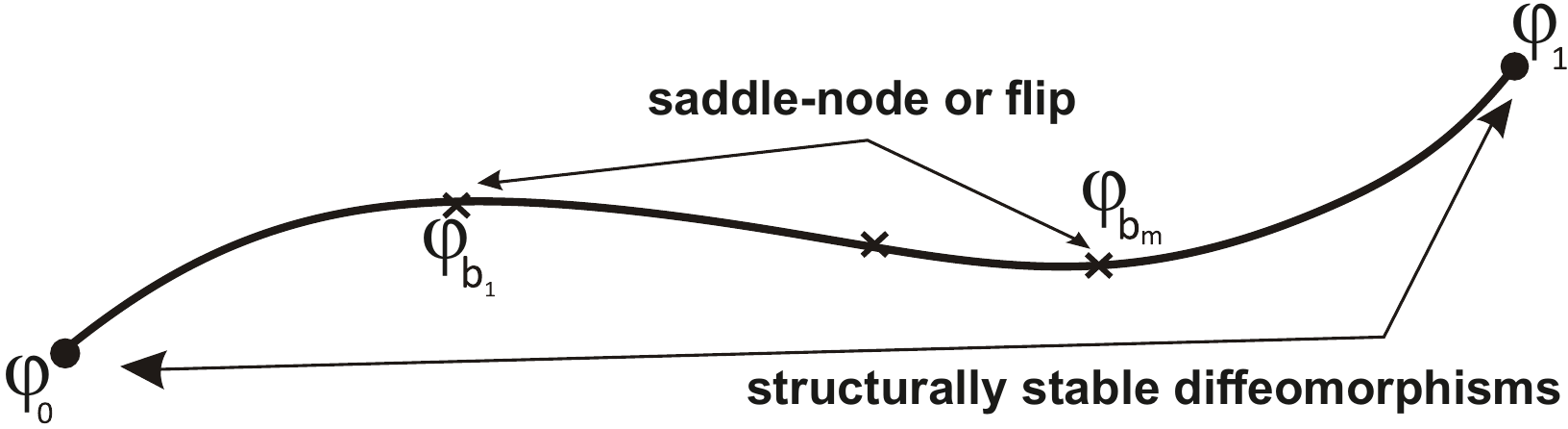}  
		\caption{\small Stable arc in the space of diffeomorphisms}	\label{arc}
	\end{figure}

	{\it Morse-Smale diffeomorphisms} (structurally stable diffeomorphisms with finite limit set) are said {\it stably connected} or {\it belong to the same  stable component} if in the space of diffeomorphisms they can be connected by a stable arc $\varphi_t$ consisting of diffeomorphisms with finite limit set. In \cite{NPT} it was found that all diffeomorphisms included to such arc are structurally stable, with the exception of a finite number of bifurcation diffeomorphisms $\varphi_{b_i},i=1,\dots,q$ (see Fig. \ref{arc}), such that:
	
	\begin{enumerate}
		
		\item the limit set of the diffeomorphism $\varphi_{b_i}$ contains a single non-hyperbolic periodic orbit which is saddle-node or flip;
		
		\item the diffeomorphism $\varphi_{b_i}$ has no cycles;
		
		\item the invariant manifolds of all periodic points of the diffeomorphism $\varphi_{b_i}$ intersect transversally;
		
		\item $\varphi_{b_i}$ has a unique nonhyperbolic periodic orbit which is the orbit of a noncritical saddle-node or of a flip which unfolds generically (see, for example, \cite{medved} for precise definitions).
	\end{enumerate}

	It is known that the simplest among Morse-Smale diffeomorphisms are the so-called {\it gradient-like diffeomorphisms}, characterized on the surface $M^2$ by the absence of intersections of stable and unstable manifolds of different saddle points. The fact that the set of such isotopic diffeomorphisms may not be stably connected, was first noted in the paper \cite{blan}. The proof of this fact was based on the construction of a special filtration for a gradient-like diffeomorphism. However, the ``sameness'' of the filtrations is only a necessary, but not sufficient condition for the existence of a stable arc between diffeomorphisms. 
	
As it turns out, a more promising approach to the stable classification of gradient-like surfaced diffeomorphisms is based on their close relation with periodic transformations. Namely, it follows from  \cite{treh} that such a diffeomorphism $f:M^2\to M^2$ is a composition
	\begin{equation}\label{comp}
		f=\phi_f\circ\xi^1_f
	\end{equation} 
of	periodic homeomorphism $\phi_f$ and a one-time shift of a Morse function gradient flow $\xi^t_f$ (see proposition \ref{mn} below).
	
	Recall that a homeomorphism $\phi :M^2 \rightarrow M^2$ is called {\it periodic} of order $m_\phi \in \mathbb{N}$ if ${\phi}^{m_\phi}=id$ and ${\phi}^j \neq id$ for any natural $j<m_\phi$. The problem of surfaced periodic homeomorphisms classification up to topological conjugacy was completely solved in  \cite{Nielsen37}, \cite{yoko}, \cite{yoko2}, \cite{yoko22}, \cite{gapa}. For a two-dimensional sphere, a complete list of periodic transformations is available in \cite{Ke} (see also \cite{cons}), for a two-dimensional torus, in \cite{BGPC}, and for orientation-preserving transformations on a surface of genus 2, in \cite{BP}.
	
	According to \cite{NoPo2024}, any isotopic to the identity gradient-like diffeomorphisms of a surface with negative Euler characteristic are connected by a stable arc.
	For surfaces with non-negative Euler characteristic it is not true. Thus, according to \cite{np2021} there are countably many stable components of gradient-like diffeomorphisms of the 2-sphere. A similar result was obtained in \cite{BaNoPo} for isotopic to the identity gradient-like diffeomorphisms of the 2-torus. 
	
	Each homeomorphism of the two-dimensional torus $\varphi:\mathbb T^2\to\mathbb T^2$ induces an action in the fundamental group uniquely determined by the matrix $\varphi_\bigstar\in GL(2,\mathbb Z)$. Moreover, the homeomorphisms $\varphi,\varphi':\mathbb T^2\to\mathbb T^2$ are diffeotopic if and only if $\varphi_\bigstar=\varphi'_\bigstar$. Thus, there are countably many isotopy classes of homeomorphisms, and each such class is uniquely determined by the matrix $A\in GL(2,\mathbb Z)$. Matrices $A,A'\in SL(2,\mathbb Z)$ are called {\it similar} over $\mathbb Z$ ($A\sim A'$) if there exists a matrix $B\in GL(2,\mathbb Z)$ such that $A'=BAB^{-1}$. For any matrix $B=\begin{pmatrix} a & c\cr
		b& d\end{pmatrix}\in GL(2,\mathbb Z)$, let $\widehat B:\mathbb T^2\to\mathbb T^2$ denote the algebraic diffeomorphism given by $$\widehat B(x,y)=(x,y)\begin{pmatrix} a & c\cr
		b& d\end{pmatrix}=(ax+by,cx+dy)\pmod 1.$$
	Let $$A_1=\begin{pmatrix} -1 & 0\cr
		0& -1\end{pmatrix},\,A_2=\begin{pmatrix} -1 & -1\cr
		1& 0\end{pmatrix},\,A_3=\begin{pmatrix} 0 & 1\cr
		-1& -1\end{pmatrix},\,A_4=\begin{pmatrix} 0 & -1\cr
		1& 1\end{pmatrix},$$
	$$A_5=\begin{pmatrix} 1 & 1\cr
		-1& 0\end{pmatrix},\,A_6=\begin{pmatrix} 0 & -1\cr
		1& 0\end{pmatrix},\,A_7=\begin{pmatrix} 0 & 1\cr
		-1& 0\end{pmatrix}.$$
	
	\begin{statement}[\cite{BGPC}, Theorem 3]\label{perA} The algebraic diffeomorphisms $\widehat A_j,\,j\in\{1,\dots,7\}$ are periodic   and any orientation-preserving periodic homeomorphism $\phi:\mathbb T^2\to\mathbb T^2$ that is not isotopic to the identity is topologically conjugate to exactly one of them. Moreover,
		\begin{equation}\label{phiA}
			\phi\,\,\hbox{is topologically conjugate to}\,\,\widehat A_j\iff \phi_\bigstar\sim A_j.
		\end{equation}
	\end{statement}
	
	Then from the relation (\ref{comp}) it follows that the set of orientation-preserving non-isotopic to the identity gradient-like diffeomorphisms of the 2-torus is divided into 7 pairwise disjoint subsets $$\mathcal G_j=\left\{f:\mathbb T^2\to\mathbb T^2:\,f_\bigstar\sim A_j\right\},\,j\in\{1,\dots,7\}.$$
	
	In the paper \cite{NoPoNo} it was established that any gradient-like diffeomorphisms
	$f,f'\in\mathcal G_1$
	are connected by a stable arc.
	
	In the present paper we consider diffeomorphisms of the set $\mathcal G_2$. For the class of diffeomorphisms under consideration we establish the following fact.
	
	\begin{lemma}\label{dyn} Any diffeomorphism $f\in \mathcal G_2$ has the  non-wandering set $\Omega_f$ with the following properties:
\begin{itemize}
			\item the set $Fix_f$ of the fixed points of $f$ consists of exactly three nodal points;
\item the set $\Omega_f\setminus Fix_f$ is non-empty, all its points have period 3, and among them there is at least one saddle orbit.
		\end{itemize}
	\end{lemma}
	
	The main result of this paper is the following theorem.
	
\begin{theorem}\label{C_i} Diffeomorphisms $f,f'\in \mathcal G_2$ are stably connected  if and only if they have the same number $i\in\{0,1,2,3\}$ of fixed sinks.
	\end{theorem}

	\section{Gradient-like diffeomorphisms of surfaces}
	
	\subsection{Basic concepts of the theory of dynamic systems}
	Let $M^n$ be a smooth closed orientable manifold and $f$ be a diffeomorphism on $M^n$.
	
	For a diffeomorphism $f$, a point $x\in M^n$ is called {\it wandering} if there exists an open neighborhood $U_x$ of $x$ such that $f^{n} (U_x)\cap U_x=\emptyset$ for all $n\in \mathbb{N}$. Otherwise, $x$ is called {\it non-wandering}.
	
	It follows immediately from the definition that every point in the neighborhood $U_x$ is wandering and, therefore, the set of wandering points is open, while the set of non-wandering points is closed.
	
	The set of all non-wandering points of a diffeomorphism $f$ is called the {\it non-wandering set} and is denoted by $\Omega_f$.
	
	The simplest examples of hyperbolic sets are, first of all, hyperbolic fixed points of a diffeomorphism, which can be classified as follows.
	Let $f:M^n\to M^n$ be a diffeomorphism and $f(p)=p$ for some $p\in M^n$. The point $p$ is {\it hyperbolic} if and only if among the eigenvalues of the Jacobian matrix $\left(\frac{\partial f}{\partial x}\right)\vert_{p}$ there are no numbers equal in absolute value to 1. If all the eigenvalues of the Jacobian matrix are less than 1 in absolute value, then $p$ is called {\it an attracting, sink point, or sink}; if all the eigenvalues are greater than 1 in absolute value, then $p$ is called {\it a repelling, source point, or source}.
	An attracting or repelling point is called {\it a nodal point}. A hyperbolic fixed point that is not a {\it nodal} is called a {\it saddle point or saddle}.
	
	If a point $p$ is a periodic point of a diffeomorphism $f$ with period $per(p)$, then, applying the previous construction to the diffeomorphism $f^{per(p)}$, we obtain a classification of hyperbolic periodic points similar to the classification of fixed hyperbolic points.
	
	The hyperbolic structure of a periodic point $p$ leads to the existence of { \it stable} $W^s_p=\{x\in M^n:\,\lim\limits_{k\to +\infty}d(f^{k\,per(p)}(x),p)\to 0\}$ and {\it unstable} $W^u_p=\{x\in M^n:\,\lim\limits_{k\to +\infty}d(f^{-k\,per(p)}(x),p)\to 0\}$ manifolds, which are smooth embeddings of $\mathbb R^{n-q_p}$ and $\mathbb R^{q_p}$, respectively. Here $q_p$ is the number of eigenvalues of the Jacobian matrix $\left(\frac{\partial f^{per(p)}}{\partial x}\right)\vert_{p}$ modulo greater than $1$, called the {\it Morse index} of the point $p$.
	
	For a hyperbolic fixed or periodic point $p$, a stable or unstable manifold is called an {\it invariant manifold} of this point, and a connected component of the set $W^u_p\setminus p$ ($W^s_p\setminus p$) is called an unstable (stable) {\it separatrix}.
	
	Recall that a diffeomorphism $f:M^n\to M^n$ is called a {\it Morse-Smale diffeomorphism} if
	
	1) the non-wandering set $\Omega_f$ consists of a finite number of hyperbolic orbits;
	
	2) the manifolds $W^s_p$, $W^u_q$ intersect transversally for any non-wandering points $p$, $q$.
	
	Let $$\Omega_q,\,q\in\{0,\dots,n\}$$ denote the subset of the nonwandering set $\Omega_f$ of the Morse-Smale diffeomorphism $f:M^n\to M^n$ consisting of points with Morse index $q$. Let $C_q=|\Omega_q|$. The symbol $\beta_q(M^n) = \beta_q$ denotes the {\it $q$-th Betti number}, that is,
	$$\beta_q(M^n) = \mathrm{rank}\, H_q(M^n, \mathbb Z).$$ By $\chi(M^n)$ we denote the Euler characteristic of $M^n$, that is,
	\[\sum\limits_{q=0}^n(-1)^q \beta_q = \chi(M^n).\]
	\begin{statement}[Lefschetz-Hopf Theorem, \cite{ShubSull_Homology75},\cite{Smale_Morse60}]\label{lefschetz} For any Morse-Smale diffeomorphism $f:M^n\to M^n$ the following relations hold:
		\begin{align}\label{beti}
			& C_0 \geq \beta_0, & C_1-C_0 \geq \beta_1 - \beta_0,\,\,\,\,\,\,\,\, & C_2 - C_1 + C_0 \geq \beta_2 -\beta_1 + \beta_0,\,\,\dots\,\,,
		\end{align}
		\begin{equation}\label{sum}
			\sum\limits_{q=0}^n(-1)^q C_q= \chi(M^n).
		\end{equation}
	\end{statement}
	
	A Morse–Smale diffeomorphism is called {\it gradient-like} if the condition $W^s_{\sigma_1} \cap W^u_{\sigma_2}\neq \emptyset $ for distinct points $\sigma_1, \sigma_2 \in \Omega_f$ implies that $dim\, W^u_{\sigma_1}< dim\, W^u_{\sigma_2}.$ In dimension $n=2$ the set of gradient-like diffeomorphisms coincides with the set of Morse–Smale diffeomorphisms whose saddle separatrices do not intersect.
	
	If $M^n$ is an orientable manifold, then a diffeomorphism $f : M^n \to M^n$ is called {\it orientation-preserving} if $f$ has a positive Jacobian at at least one point, otherwise the diffeomorphism is called {\it orientation-changing}.
	
	\subsection{Attractors and Repellers}	
	Recall that a compact $f$-invariant set $A\subset M^n$ is called an {\it attractor} of a diffeomorphism $f:M^n\to M^n$ if it has a compact neighborhood $U_{A}$ such that $f(U_{A})\subset {\rm int}\,(U_{A})$ and $A = \underset{k\geqslant 0}{\bigcap} f^{k}(U_{A}).$ The neighborhood $U_{A}$ is called {\it capturing}. {\it Repeller} is defined as an attractor for the diffeomorphism $f^{-1}$.
	
\begin{statement}[\cite{GrZhuMePo}, Theorem 1.1]\label{<n-2} Let $\Sigma\subset\Omega_f$ be an $f$-invariant set for a Morse-Smale diffeomorphism $f:M^n\to M^n$ such that the union of the unstable manifolds of all points from $\Sigma$ is closed. Then
	\begin{itemize}
\item $A_\Sigma=\bigcup\limits_{p\in\Sigma}W^u_p$ is an attractor of  $f$,
\item $R_\Sigma=W^s_{\Omega_f\setminus\Sigma}$ is a repeller of $f$,
\item $A_\Sigma$ is connected if $\dim\,R_\Sigma\leqslant n-2$.	
\end{itemize}	  
	\end{statement}
	
The attractor $A_\Sigma$ and the repeller $R_\Sigma$ are called {\it dual}.

\begin{statement}[Lemmas 3.2, 8.1  \cite{np2021}]\label{sad-noo} Let $f:M^2\to M^2$ be a gradient-like diffeomorphism and  there are pairwise disjoint 2-disks
			$$D,f(D),\dots,f^{m-1}(D),\,m\in\mathbb N$$ such that $f^m(D)\subset {\rm int}\, D$, and $f^m$ has a fixed sink in $D$. Then there exists a stable arc connecting $f$ with a a gradient-like diffeomorphism  $\tilde f$, which coincides with $f$ out of  the discs and $\tilde f^m$ has a unique non-wandering point in $D$ -- the fixed sink.
	\end{statement}	
	
\begin{statement}[Lemma 8.5  \cite{np2021}]\label{ann} Let $f:M^2\to M^2$ be a gradient-like diffeomorphism and there is an closed annulus $Q$ such that $f(Q)\subset {\rm int}\, Q$. Then there exists a stable arc connecting $f$ with a a gradient-like diffeomorphism  $\tilde f$, which coincides with $f$ out of  $Q$ and $Q$ is a trapping neighborhood for an attractor $\tilde A$ of $\tilde f$ such that $\tilde f|_{\tilde A}$ is topologically conjugate to a rough transformation of the circle with a unique sink and a unique source orbits.
	\end{statement}	
		
	\subsection{Relationship with periodic homeomorphisms}
	
	Let $M^2$ be a closed connected orientable surface.
	For any homeomorphism $\phi:M^2 \rightarrow M^2$, denote by $Per_\phi$ the set of its periodic points, by $P_\phi$ the set of periods of periodic points, and by $Fix_\phi$ the set of its fixed points.  
	
	Recall that a homeomorphism $\phi :M^2 \rightarrow M^2$ is called {\it periodic} of order $m_\phi \in \mathbb{N}$ if ${\phi}^{m_\phi}=id$ and ${\phi}^j \neq id$ for any natural $j<m_\phi$. It follows from the definition of a periodic homeomorphism $\phi$ that $P_\phi$ is finite and consists of divisors of $m_\phi$. Set $\tilde P_\phi=P_\phi \backslash\{m_\phi\}$. For any $l\in\tilde P_\phi$, denote by $B_\phi^l$ the set of points of period $l$. Set $B_\phi=\bigsqcup\limits_{l\in \tilde P_\phi} B_\phi^l$.

	\begin{statement}[\cite{Nielsen37}, \cite{yoko}]\label{data}
		Let $\phi :M^2 \rightarrow M^2$ be an orientation-preserving periodic homeomorphism of order $m_\phi$. Then
		$B_\phi^l$ is finite for any $l\in \tilde P_{\phi}$.
	\end{statement}	
	
	\begin{statement}[\cite{treh}]\label{mn} Let $f:M^2 \rightarrow M^2$ be an orientation-preserving gradient-like diffeomorphism. Then there exist an orientation-preserving periodic homeomorphism $\phi_f:M^2 \rightarrow M^2$ of order $m_{f}$ and a flow $\xi^t_f:M^2 \rightarrow M^2$ equivalent to the gradient flow of a generic Morse function such that
		\begin{enumerate}
			\item[1)] $f=\phi_f\circ\xi^1_f,\, \phi_f\circ\xi^t_f=\xi^t_f\circ\phi_f$;
			\item[2)] $B_{\phi_f}\subset\Omega_f=Fix_{\xi^1_f}$ and the invariant manifolds of periodic points of the diffeomorphism $f$ coincide with the invariant manifolds of fixed points of the flow $\xi^t_f$;
			\item[3)] $f|_{\Omega_f}=\phi_f|_{\Omega_f}$ and $\tilde P_{\phi_f} \subset P_f \subset P_{\phi_f}$;
			\item[4)] $m_f$ is the smallest of the natural numbers $m$ such that the set $\Omega_f$ consists of fixed points of the diffeomorphism $f^m$ and all saddle points have positive orientation type;
			\item[5)] the period of any saddle separatrix is equal to $m_f$.
		\end{enumerate}
	\end{statement}	
	
	\section{Dynamics of diffeomorphisms of the set $\mathcal G_2$}
	In this section we prove Lemma \ref{dyn}. To do this, recall that we set $A_2=\begin{pmatrix} -1 & -1\cr
		1& 0\end{pmatrix}$ and introduced the set of gradient-like diffeomorphisms $$\mathcal G_2=\left\{f:\mathbb T^2\to\mathbb T^2:\,f_\bigstar\sim A_2\right\}.$$
	Consider the following subset of this set: $$G_2=\left\{f:\mathbb T^2\to\mathbb T^2:\,f_\bigstar=A_2\right\}.$$
	
	\begin{lemma}\label{sopr} Any diffeomorphism $f\in\mathcal G_2$ is topologically conjugate to some diffeomorphism $g\in G_2$.
	\end{lemma}
	\begin{proof} Let $f_\bigstar=A'_2$. Then $A_2\sim A'_2$. Hence there exists a matrix $B\in GL(2,\mathbb Z)$ such that $A_2=BA'_2B^{-1}$. Consider an algebraic diffeomorphism $\widehat B:\mathbb T^2\to\mathbb T^2$ and set $g = \widehat B f \widehat B^{-1}$. Then the diffeomorphisms $f$ and $g$ are topologically conjugate, and $g_\bigstar=A_2$.
	\end{proof}
	
	\begin{statement}[\cite{BGPC}, Theorem 1]\label{A2} The periodic diffeomorphism $\widehat A_2$ has the following properties:
		\begin{itemize}
			\item $m_{\widehat A_2}=3$;
			\item the set $B_{\widehat A_2}$ consists of three fixed points;
			\item in the neighborhood of each fixed point the diffeomorphism $\widehat A_2$ is topologically conjugate to the rotation of the plane $z\mapsto ze^{i\frac{2\pi}{3}}$.
		\end{itemize}
	\end{statement}

	It follows directly from Lemma \ref{sopr} that the proof of Lemma \ref{dyn} is reduced to the proof of the following statement.
	
	\begin{lemma}\label{dyng} Any diffeomorphism $f\in G_2$  has the  non-wandering set $\Omega_f$ with the following properties:
\begin{itemize}
			\item the set $Fix_f$ of the fixed points of $f$ consists of exactly three nodal points;
\item the set $\Omega_f\setminus Fix_f$ is non-empty, all its points have period 3, and among them there is at least one saddle orbit.
\end{itemize}
	\end{lemma}
	\begin{proof} By Statement \ref{mn}, the diffeomorphism $f$ can be represented as a composition $f=\phi_f\circ\xi^1_f$ of a periodic homeomorphism $\phi_f$ and a shift by unit time of the flow. It follows that $f$ and $\phi_f$ are isotopic and, consequently, $\phi_{f\bigstar}=A_2$. It follows from the criterion (\ref{phiA}) of Statement \ref{perA} that the homeomorphism $\phi_f$ is topologically conjugate to the diffeomorphism $\widehat A_2$. Then, by Statement \ref{A2}, the period of the homeomorphism $\phi_f$ is equal to 3 and the set of its points of smaller period $B_{\phi_f}$ consists of exactly three fixed points.
		
From items 2) and 3) of Statement \ref{mn} it follows that the nonwandering set $\Omega_f$ contains exactly three fixed points, all its other points have period $3$. From item 5) of Statement \ref{mn} it follows that any saddle point of the diffeomorphism $f$ has period $3$, and hence the set $Fix_f$ consists of three nodal points. Since the Euler characteristic of a two-dimensional torus is 0, it follows from equality \eqref{sum} of the Lefschetz-Hopf theorem (see Statement \ref{lefschetz}) that the set $\Omega_f\setminus Fix_f$ is nonempty and contains at least one saddle point. From item 3) of Statement \ref{mn} it follows that all points of the set $\Omega_f\setminus Fix_f$ have period 3. 
	\end{proof}
	
	For $i\in\{0,\dots,3\}$ we set $$G_{2,i}=\{g\in G_2:\, Fix_g\,\,\hbox{contains exactly}\,\,i\,\,\hbox{fixed sinks}\}.$$
	
	A diffeomorphism $g_i\in G_{2,i}$ is called {\it the simplest} if $$|\Omega_{g_i}|=\min_{g\in G_{2,i}}|\Omega_g|.$$

	\section{Dynamics of the simplest diffeomorphisms}
	In this section we describe the dynamic properties of the simplest diffeomorphisms $g_i\in G_{2,i},\,i\in\{0,1,2,3\}$. In what follows they will play the role of representatives of stable components in the set $G_2$.   
	
	\begin{lemma}\label{g1} Any diffeomorphism $g_1\in G_{2,1}$ has the following properties (see Fig. \ref{picg1}):
		\begin{enumerate}
			\item the set $\Omega_{g_1}$ consists of two fixed sources $\alpha_1,\alpha_2$, one fixed sink $\omega$, and one saddle orbit $\sigma,f(\sigma),f^2(\sigma)$;
			\item the set $K=W^u_\sigma\cup\omega$ is a non-contractible knot on the torus $\mathbb T^2$;
			\item knots $K,g_1(K),g_1^2(K)$ have homotopy types $\pm\langle 1,0\rangle$, $\mp\langle 1,1\rangle$, $\pm\langle 0,1\rangle$.
		\end{enumerate}
		Any two diffeomorphisms of the set $G_{2,1}$ are topologically conjugate.
	\end{lemma}
	\begin{proof} We prove each point of the lemma one after another.
		
		1. By the Lemma \ref{dyng} and the definition of the diffeomorphism $g_1$, the set $\Omega_{g_1}$ contains two fixed sources $\alpha_1,\alpha_2$, one fixed sink $\omega$, and at least one saddle orbit $\sigma,g_1(\sigma),g_1^2(\sigma)$ of period 3. The minimality condition and formula \eqref{sum} indicate that this saddle orbit is unique in the set $\Omega_{g_1}$.
		
		2. Since the set $\Omega_{g_1}$ contains a unique sink $\omega$, then ${\rm cl}\,W^u_\sigma\setminus W^u_\sigma=\omega$. Therefore, the set $K=W^u_\sigma\cup\omega$ is a knot on the torus $\mathbb T^2$. Let us show that it is not contractible. Assuming the contrary, we obtain that the knot $K$ bounds the 2-disk $D\subset\mathbb T^2$. By item 5) of Statement \ref{mn}, each saddle separatrix has period 3, which implies that the disks $D,g_1(D),g_1^2(D)$ intersect at a unique point $\omega$. On the other hand, the disk $D$ contains a stable separatrix of the saddle $\sigma$, and therefore contains some source in its closure. Then the period of such a source is at least three, which contradicts the assumption about the structure of the set $\Omega_{g_1}$.
		
		Thus, $K$ is not contractible knot.
		
		3. Let $K$ be a knot of homotopy type $\langle a,b\rangle\neq\langle 0,0\rangle$. Since $g_{1\bigstar}=A_2$, similarly to Lemma \ref{g1}, it is possible to prove that  the knots $K,g_1(K),g_1^2(K)$ have the following homotopy types:
		\begin{itemize}
			\item $\langle K\rangle=\langle a,b\rangle$;
			\item $\langle g_1(K)\rangle=\langle a,b\rangle A_2=\langle b-a,-a\rangle$;
			\item $\langle g_1^2(K)\rangle=\langle b-a,-a\rangle A_2=\langle -b,a-b\rangle$.
		\end{itemize}
		By virtue of the proved point 2, $K,g_1(K),g_1^2(K)$ are not contractible and have a unique intersection point $\omega$. According to \cite{Rol}, the determinants of the matrices of all matrices $\begin{pmatrix} a & b\cr
			b-a& -a\end{pmatrix},$ $\begin{pmatrix} b-a & -a\cr
			-b& a-b\end{pmatrix},$ $\begin{pmatrix} -b & a-b\cr
			a& b\end{pmatrix}$ in this case must be equal to $\pm 1$. By direct calculation we see that this condition is equivalent to the equality $$-a^2 + ab -b^2 = \pm1.$$
		
Let's consider the case when $-a^2 + ab -b^2 = 1$. By representing the equation as a quadratic one in the variable $a$, we get that the discriminant of this equation is $-3b^2 - 4$. It is always negative, which means that in this case no pair of real numbers $\langle a, b\rangle$ can be its solution. 
		
	In the case of $-a^2 + ab -b^2 = -1$, the discriminant of the equation, as a square equation with respect to $a$, is equal to $4 - 3b^2$. In this case, the condition $4 - 3b^2\geqslant 0$
		is satisfied by only three integer values of $b$: $0$, $1$, and $-1$.
		Taking into account that, by virtue of point 2, $\langle a,b\rangle\neq\langle0,0\rangle$, we obtain the following possible homotopy types of the knot $\langle a,b\rangle$: $\pm\langle 1,0\rangle$, $\mp\langle 1,1\rangle$, $\pm\langle 0,1\rangle$.
		
To depict the phase portrait of the diffeomorphism $g_1$, we arrange on the torus $\mathbb T^2$ knots $K,g_1(K),g_1^2(K)$ in accordance with their homotopy types $\langle 1,0\rangle$, $\langle -1,-1\rangle$, $\langle 0,1\rangle$ (see Fig. \ref{picg1}).
		\begin{figure}[h!]		\centerline{\includegraphics
				[width=7 cm]{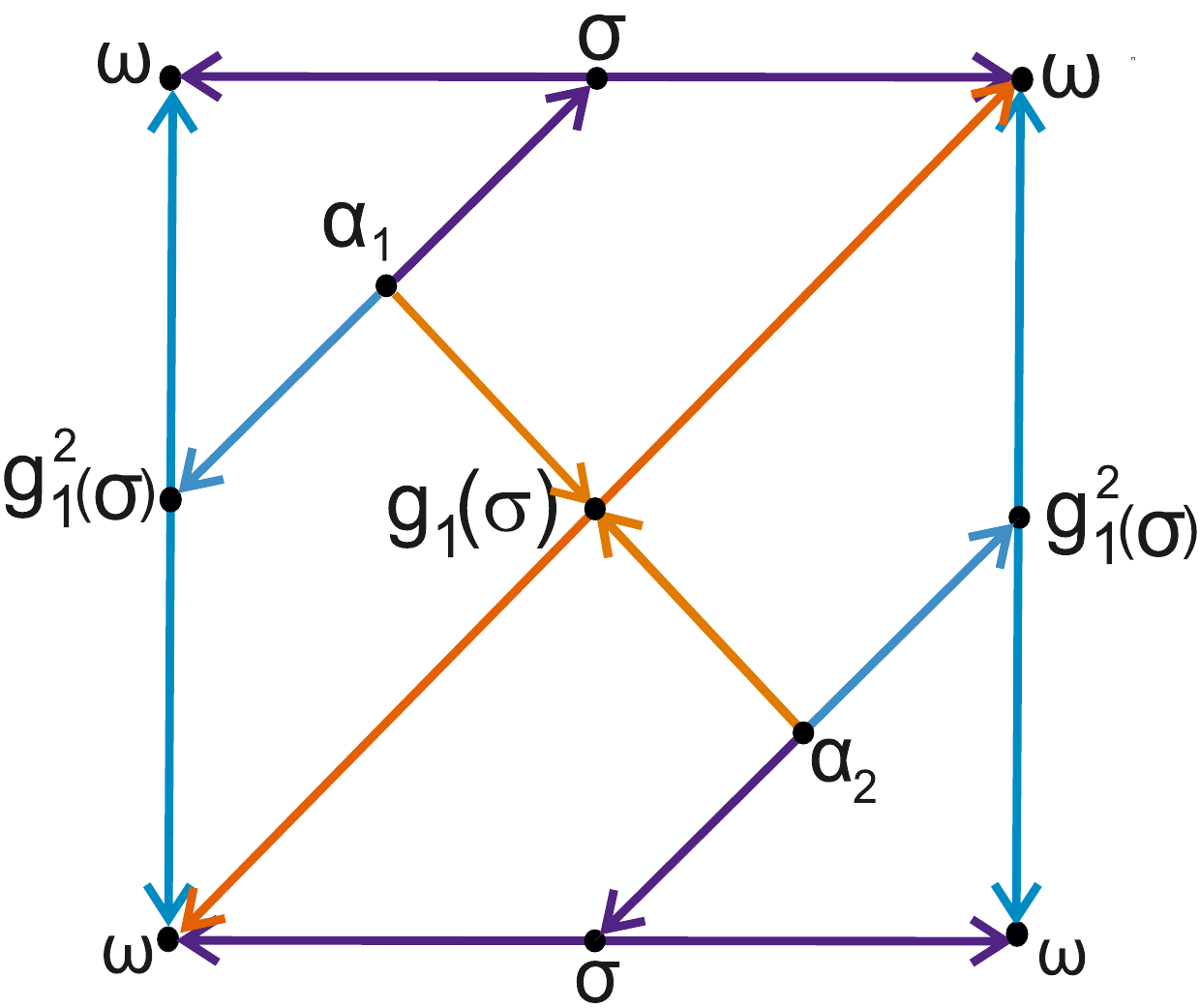}}		\caption{\small Phase portrait of the diffeomorphism $g_1$}\label{picg1}	
		\end{figure}	
		
		It remains to show that any two diffeomorphisms of the set $G_{2,1}$ are topologically conjugate. To do this, we use a three-color graph, which is a complete invariant of gradient-like diffeomorphisms (see \cite{treh}). To construct the graph, we color all stable (unstable) saddle separatrices blue (red). In each region complementary to the closure of the saddle separatrices, we choose one invariant curve and color it green. As a result, the torus $\mathbb T^2$ is divided by colored curves into triangular regions. We assign a vertex to each such region and connect two vertices with an edge of the corresponding color if the regions have a boundary of this color. The resulting graph $\Gamma_{g_1}$ is called the {\it three-color graph} of the diffeomorphism $g_1$ (see Fig. \ref{gr1}). 
		\begin{figure}[h!]		\centerline{\includegraphics
				[width=12 cm]{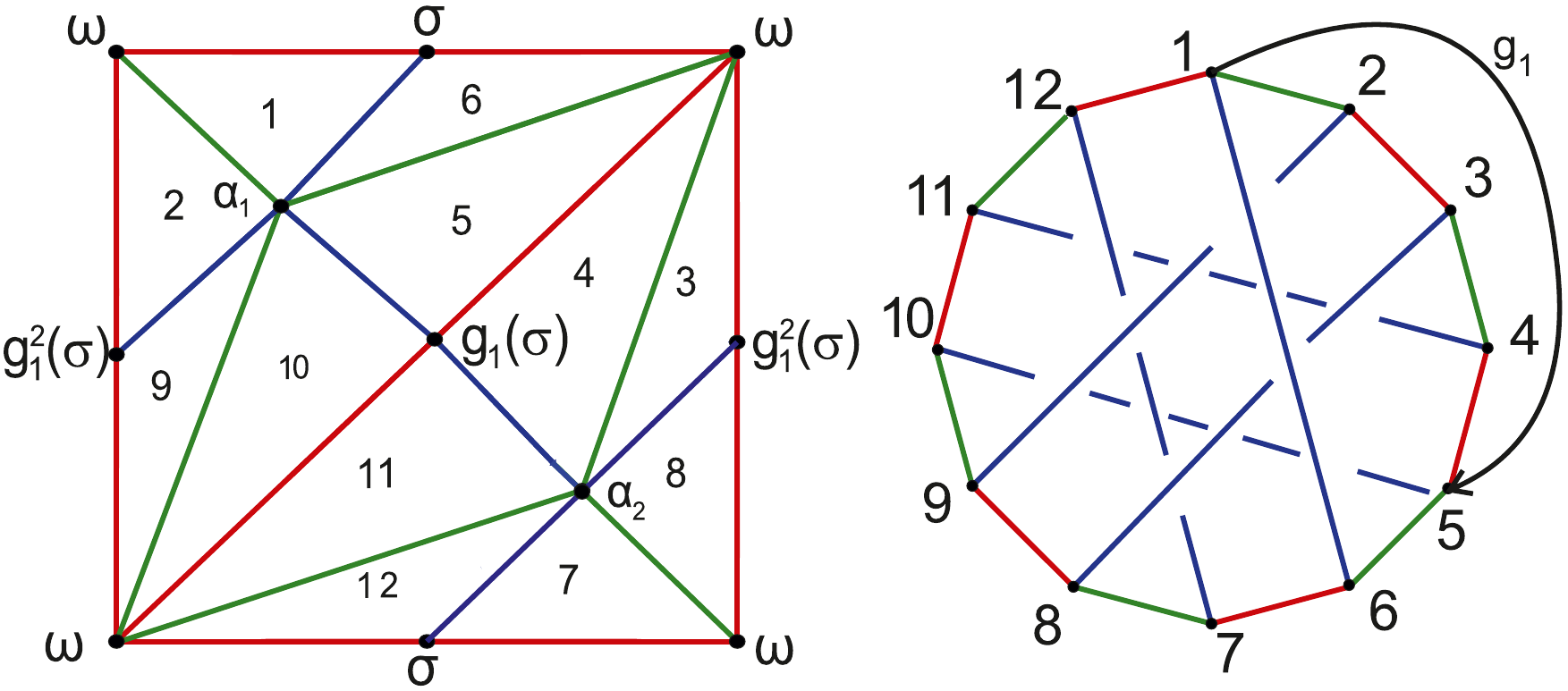}}		\caption{\small Three-color  graph of diffeomorphism $g_1$}\label{gr1}	
		\end{figure}	
		
		By construction, all vertices of the graph $\Gamma_{g_1}$ form a unique red-green cycle, and the diffeomorphism $g_1$ induces a permutation of the graph vertices consisting of rotating this cycle by an angle $\frac{2\pi}{3}$. Two three-color graphs {\it are equivalent} if there exists an isomorphism between them that preserves the color of edges and conjugates permutations. The equivalence class of a graph is a complete invariant of the topological conjugacy of the corresponding gradient-like diffeomorphism. Moreover, two such graphs are obviously equivalent, which completes the proof.
	\end{proof}
	
	\begin{lemma}\label{g0} Any diffeomorphism $g_0\in G_{2,0}$ has the following properties (see Fig. \ref{picg0}):
		\begin{enumerate}
			\item the set $\Omega_{g_0}$ consists of three fixed sources $\alpha_1,\alpha_2,\alpha_3$, one sink orbit $\omega,g_0(\omega),g_0^2(\omega)$ and two saddle orbits $\sigma_1,g_0(\sigma_1),g_0^2(\sigma_1)$, $\sigma_2,g_0(\sigma_2),g_0^2(\sigma_2)$;
			\item the set $K=\omega\cup W^u_{\sigma_1}\cup W^u_{\sigma_2}\cup g_0(\omega)$ is non-contractible knot on the torus $\mathbb T^2$;
\item knots $K,g_0(K),g_0^2(K)$ have homotopy types $\pm\langle 1,0\rangle$, $\mp\langle 1,1\rangle$, $\pm\langle 0,1\rangle$.
		\end{enumerate}
		Any two diffeomorphisms of the set $G_{2,0}$ are topologically conjugate.
	\end{lemma}
	\begin{proof}
		Let us prove each point of the lemma sequentially.

1. By Lemma \ref{dyng} and the definition of the diffeomorphism $g_0$, the set $\Omega_{g_0}$ contains three fixed sources $\alpha_1,\alpha_2, \alpha_3$, at least one sink orbit of period 3, and at least one saddle orbit of period 3. The minimality condition and formula \eqref{sum} indicate that the sink orbit $\omega,g_0(\omega),g_0^2(\omega)$ is unique in the set $\Omega_{g_0}$, and there are two saddle orbits: $\sigma_1,g_0(\sigma_1),g_0^2(\sigma_1)$, $\sigma_2,g_0(\sigma_2),g_0^2(\sigma_2)$.
		
2. Let us show that unstable saddle separatrices of saddles $\sigma_1,g_0(\sigma_1),g_0^2(\sigma_1)$, $\sigma_2,g_0(\sigma_2),g_0^2(\sigma_2)$ go to different points of the sink orbit. Let us assume the opposite. Let, for definiteness, two separatrices of the saddle $\sigma_1$ be connected to the sink $\omega$. Then consider the knot $\Gamma = W^u_\sigma\cup\omega$, which will be non-contractible on the torus. If it was contractible, then there would be a source orbit of the period $3$, that contradicts item 1. On the other side  $\Gamma$ and $g_0(\Gamma)$ do not intersect and, hence, have the same homotopy type. The only homotopy type that matrix $A_2$ preserves is $(0,0)$, but $\Gamma$ is non-contractible, it is a contradiction. 

Thus, the unstable saddle separatrices of saddles connect to different points of the sink orbit. Then the set $K=\omega\cup W^u_{\sigma_1}\cup W^u_{\sigma_2}\cup g_0(\omega)$ is a knot on the torus. This knot is non-contractible as above.
		
3. Let $K$ be a knot of homotopy type $\langle a,b\rangle\neq\langle 0,0\rangle$. Since $g_{0\bigstar}=A_2$, similarly to item 3 of Lemma \ref{g1}, it is possible to prove that  the knots $K,g_0(K),g_0^2(K)$ have the following homotopy types: $\pm\langle 1,0\rangle$, $\mp\langle 1,1\rangle$, $\pm\langle 0,1\rangle$.
		
To depict the phase portrait of the diffeomorphism $g_0$, we arrange on the torus $\mathbb T^2$ knots $K,g_0(K),g_0^2(K)$ in accordance with their homotopy types $\langle 1,0\rangle$, $\langle -1,-1\rangle$, $\langle 0,1\rangle$ (see Fig. \ref{picg0}).
			\begin{figure}[h!]		\centerline{\includegraphics
				[width=7 cm]{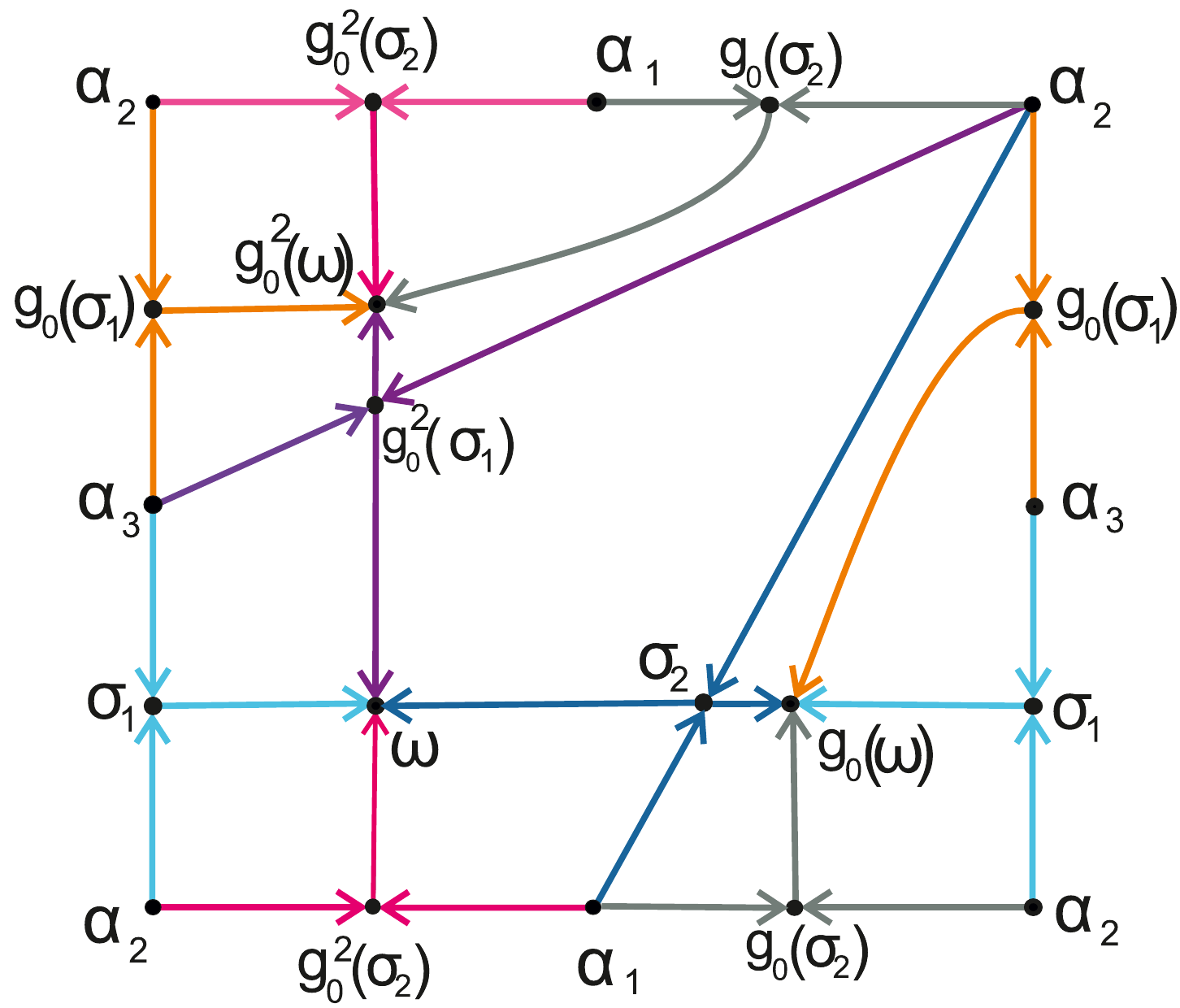}}		\caption{\small Phase portrait of the diffeomorphism $g_0$}\label{picg0}	
		\end{figure}	
		
As in the case of the diffeomorphism $g_1$, the topological conjugacy of any two diffeomorphisms of the set $G_{2,0}$ is proved via a three-color graph (see Fig. \ref{gr0}).	
	\begin{figure}[h!]		\centerline{\includegraphics
			[width=12 cm]{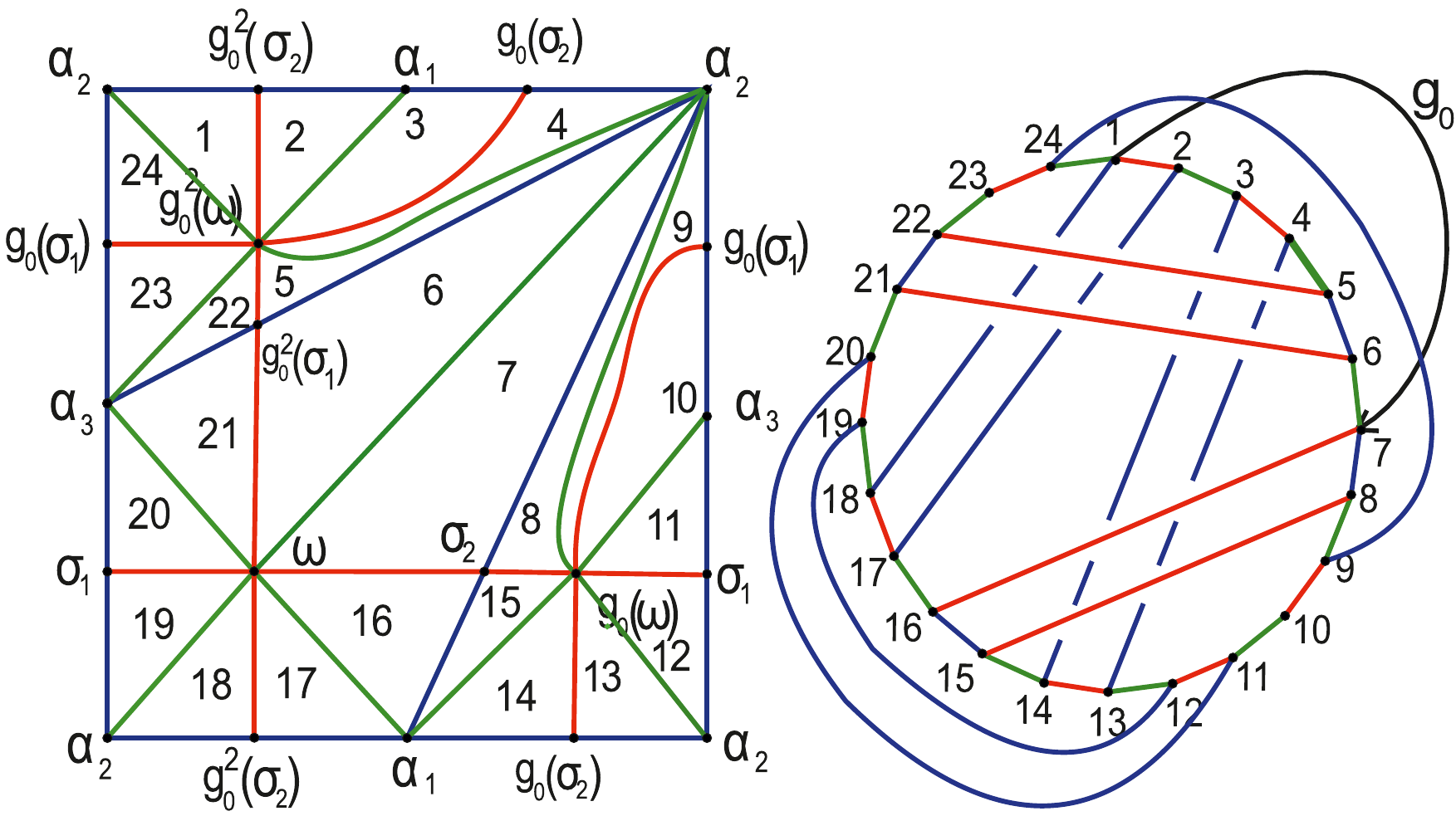}}		\caption{\small Three-color graph of diffeomorphism $g_0$}\label{gr0}	
	\end{figure}	
\end{proof}
	
\begin{figure}[h!]		\centerline{\includegraphics
		[width=7 cm]{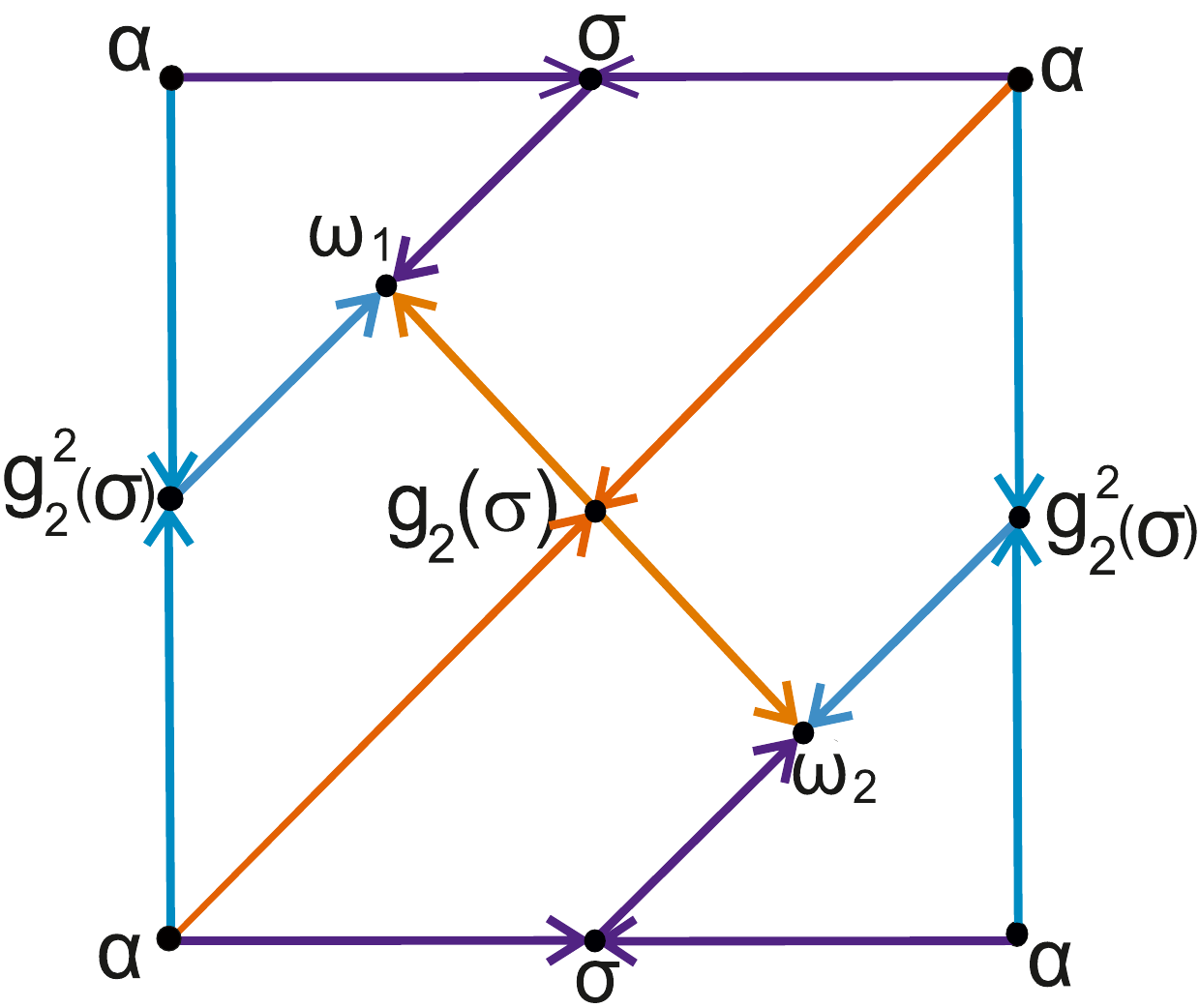}}		\caption{\small Phase portrait of the diffeomorphism $g_2$}\label{picg2}	
\end{figure}

\begin{figure}[h!]		\centerline{\includegraphics
		[width=7 cm]{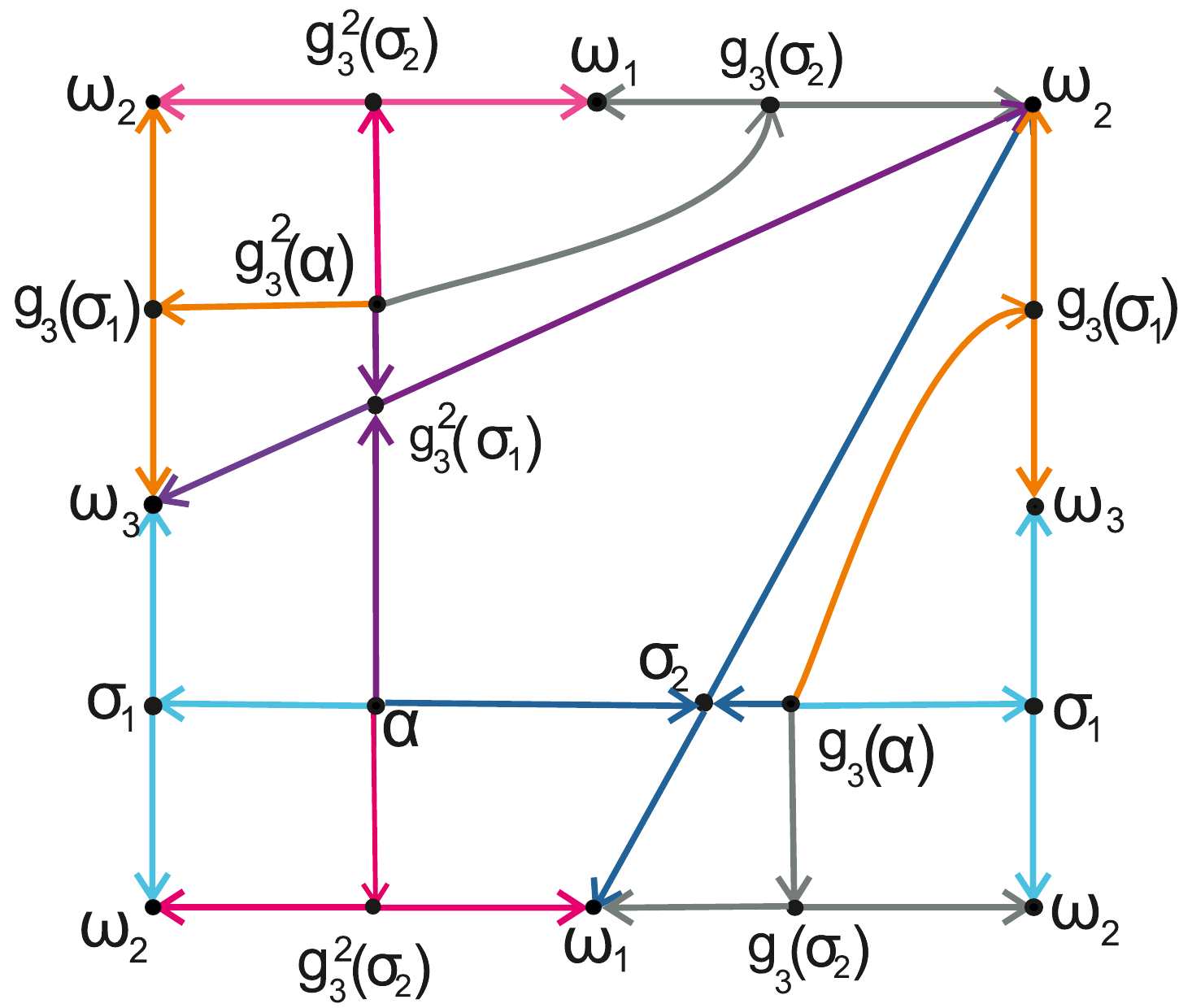}}		\caption{\small Phase portrait of the diffeomorphism $g_3$}\label{picg3}	
\end{figure}

It follows directly from the definition of the simplest diffeomorphisms that the phase portrait of every simplest  diffeomorphism $g_2\in G_{2,2}$ (see Fig. \ref{picg2}) is obtained from the phase portrait of the simplest diffeomorphism $g_1\in G_{2,1}$ by formally replacing sinks with sources and vice versa, as well as by changing the orientation of the saddle separatrices. Similar rule works for the simplest diffeomorphisms $g_3\in G_{2,3}$ (see Fig. \ref{picg3}) and $g_0\in G_{2,0}$.

	\section{Stable components of class $G_2$}
	In this section we prove Theorem \ref{C_i}. The proof is conceptually divided into establishing the following facts.
	\begin{enumerate}
		\item Diffeomorphisms of classes $G_{2,i},G_{2,j},\,i\neq j$ are not connected by a stable arc (Lemma \ref{step1} below).
		\item Any diffeomorphism $g\in G_{2,i}$ is connected by a stable arc with some the simplest diffeomorphism $g_i\in G_{2,i}$ (Lemma \ref{step3} below).
		\item Any the simplest diffeomorphisms of class $G_{2,i}$ are connected by an arc without bifurcations (Lemma \ref{step2} below).
	\end{enumerate}
	
	To prove Lemma \ref{step1} we need the following observation.
	
Let $\omega$ be a fixed sink point with the positive type of the orientation for a gradient-like diffeomorphism $f=\phi_f\circ\xi^1_f:M^2\to M^2$ and $c_\omega$ be a section to the trajectories of the flow $\xi^t_f$ in  $W^s_\omega\setminus\omega$. Since $\phi_f$ implements the equivalence of the flow $\xi^t_f$ with itself, it induces a homeomorphism $$\phi_\omega:c_\omega\to c_\omega,$$ which maps the intersection point of the trajectory $\mathcal O$ of the flow $\xi^t_f$ with the section $c_\omega$ to the intersection point of the trajectory $\phi_f(\mathcal O)$ with one. Denote by $\rho_\omega$ the rotation number of the homeomorphism $\phi_\omega$ and call it {\it the rotation number of the sink $\omega$}. 
	
\begin{lemma}\label{rotw} Let $\omega$ be a fixed sink point with the positive type of the orientation for a gradient-like diffeomorphism $f=\phi_f\circ\xi^1_f:M^2\to M^2$ and  $f'=\phi_{f'}\circ\xi^1_{f'}:M'^2\to M'^2$ be gradient-like diffeomorphisms topologically conjugate to $f$ via a homeomorphism $h:M^2\to M'^2$. Then
\begin{equation}\label{rh}
	\rho_{h(\omega)}\equiv\pm\rho_{\omega}\pmod 1.
		\end{equation}
			\end{lemma}
	\begin{proof} If the diffeomorphism $f$ has no saddle points then, by \cite[Theorem 2.5]{grin}, $f$ is a ``source-sink'' diffeomorphism. In this case, by Statement \ref{mn}, $\phi_f=id$ and, hence,  $\rho_{\omega}=0$. As $f'$ is topologically conjugate to $f$ by means $h$ then $\phi_{f'}=id$ and $\rho_{h(\omega)}=0$.
		
If a diffeomorphism $f$ has at least one saddle point, then by \cite[Corollary 2.2]{grin} there exists at least one saddle point whose unstable separatrix $l$ belongs to the basin $W^s_\omega$. Denote by $L_\omega$ the union of all such separatrices. By Statement \ref{mn}, each of these separatrices $l$ is a trajectory of the flow $\xi^t_f$, that implies that $l\cap c_\omega$ consists of a unique point.  
		
Since $h$ conjugates $f$ with $f'$ then $\omega'$ is  a fixed sink point with the positive type of the orientation for $f'$ and $L_{\omega'}=h(L_\omega)$ is the union of all unstable saddle separatrices $l'$ belonging to the basin $W^s_{\omega'}$. Moreover, each of the intersections $l'\cap c_{\omega'},\,l'\cap h(c_{\omega})$ consists of one point (see Fig. \ref{rov}). 
		\begin{figure}[h!]
			\centerline{\includegraphics [width=8 cm]{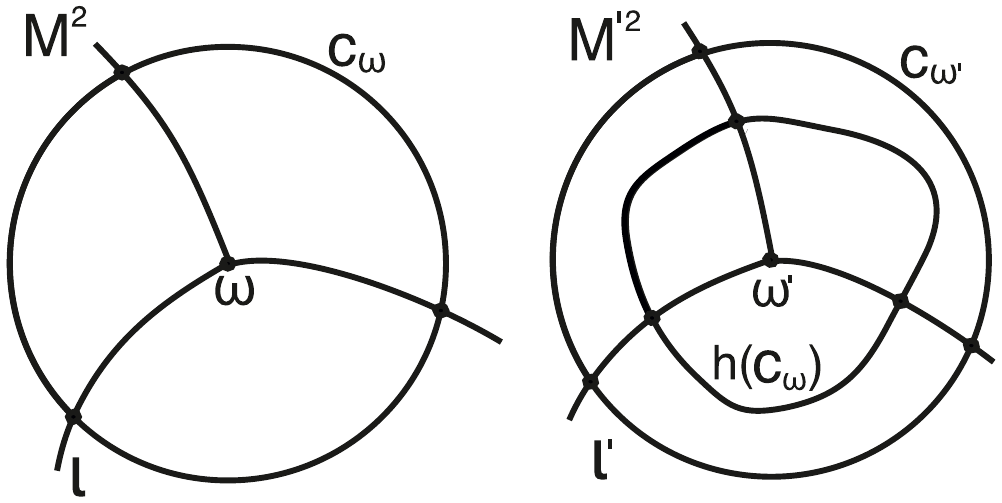}}
			\caption{\small Illustration to  Lemma \ref{rotw}}\label{rov}
		\end{figure}  
		
Let $H:c_\omega\to c_{\omega'}$ be a homeomorphism such that for any connected component $\gamma$ of the set $c_\omega\setminus L_\omega$, the arcs $h(\gamma)$ and $H(\gamma)$ belong to the same connected component of the set $W^s_{\omega'}\setminus (L_{\omega'}\cup\omega')$. Then 
		\begin{equation}\label{Hp}
			H\phi_\omega|_{c_\omega\cap L_\omega}=\phi_{\omega'}H|_{c_\omega\cap L_\omega}.
		\end{equation} 
		Since the rotation number of a homeomorphism of a circle is determined by one of its orbits (see, for example, \cite[Proposition 11.1.1]{Kat}), then, by virtue of \cite[Proposition 11.1.3]{Kat}, the relation \eqref{Hp} implies the equality \eqref{rh}. 
	\end{proof}

A fact similar to Statement \ref{rotw} takes place for  a fixed source point with the positive type of the orientation.

\begin{lemma}\label{step1} For any stable arc $f_t,\,t\in[0,1]$ such that $f_0\in G_{2,i}$, it is true that $f_1\in G_{2,i}$.
	\end{lemma}
	\begin{proof} To prove this, by Lemma \ref{dyn}, it suffices to show that passing along a stable arc does not destroy the set $Fix_{f_0}$. Indeed, to destroy $Fix_{f_0}$ by saddle-node or flip bifurcation we need an $f_0$-invariant curve (central manifold) passing trough the point of this set.  However, by Lemma \ref{rotw}, all points of the set $Fix_{f_0}$ have the rotation number $\frac13$ and, hence, such a curve does not exist. As any bifurcation on the stable arc preserves dynamics (up to topological conjugacy) out of a neighborhood of the bifurcated points, then a diffeomorphism $f_t$ in a neighborhood of ${Fix_{f_t}}$ is topologically conjugate to the diffeomorphism $f_0$ in a neighborhood of ${Fix_{f_0}}$ for every $t\in[0,1]$. 
\end{proof}
	
\begin{lemma}\label{step3} Any diffeomorphism $g\in G_{2,i}$ is connected by a stable arc to some the simplest diffeomorphism.
\end{lemma}
\begin{proof} We prove the lemma for the cases $i=0,1$, for other cases the proof is reduce to these one by the transition to the inverse maps. 

{\bf $i=1$.} Let $g\in G_{2,1}$. According to the Lemma \ref{dyng}, the set $Fix_g$ consists of one sink $\omega$ and two sources $\alpha_1,\alpha_2$. Let's consider the set $\Sigma_{\omega} = \{\sigma: \omega \in {\rm cl}\,W^u_{\sigma}\}$. This set is not empty (see, for example, \cite[Corollary 2.2]{grin}). Let $$\Gamma_{\omega}=\bigcup\limits_{\sigma\in \Sigma_\omega}{\rm cl}\,W^u_{\sigma}.$$ For one-dimensional complex $\Gamma_{\omega}$ there are two  possible cases: 1) $\Gamma_{\omega}$ contains a non-contractible loop, 2) other case. Let's consider these cases separately. 

1) Denote by $K_\omega$ a non-contractible loop in $\Gamma_{\omega}$. It can be of two types: either $K_\omega={\rm cl}\,W^u_\sigma$ for a saddle $\sigma\in\Sigma_\omega$, or $K_\omega={\rm cl}\,(W^u_{\sigma_1}\cup W^u_{\sigma_2})$ for saddles $\sigma_1,\sigma_2\in\Sigma_\omega$ (see Fig. \ref{Aw}).  
\begin{figure}[h!]
			\centerline{\includegraphics [width=7 cm]{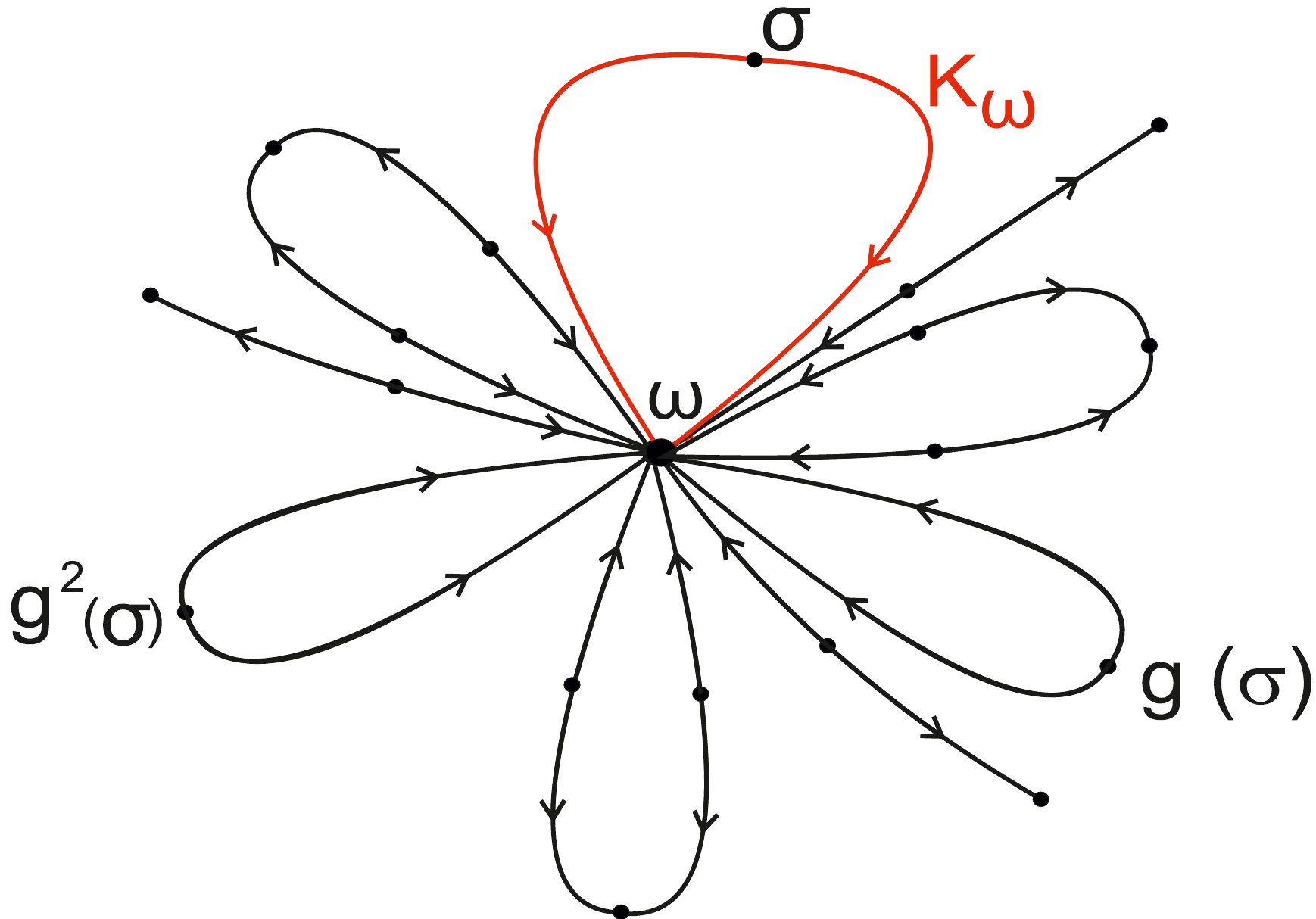}}
			\caption{\small Complex $\Gamma_\omega$ with non-contractible loop}\label{Aw}
		\end{figure} 
Up to saddle-node bifurcation, we can assume that $K_\omega={\rm cl}\,W^u_\sigma$. As $g_{\bigstar}=A_2$, similarly to Lemma \ref{g1}, it is possible to prove that  the knots $K_\omega,g(K_\omega),g^2(K_\omega)$ have the following homotopy types: $\pm\langle 1,0\rangle$, $\mp\langle 1,1\rangle$, $\pm\langle 0,1\rangle$. Moreover, by Statement \ref{<n-2}, $$A_\omega=K_\omega\cup g(K_\omega)\cup g^2(K_\omega)$$ is an attractor of the diffeomorphism $g$. Then there is its trapping  neighborhood $U_{A_{\omega}}$ such that  $\mathbb T^2 \setminus{\rm int}\,U_{A_{\omega}}$ is the disjoint union of two 2-discs $D_1 \sqcup D_2$ and $g^{-1}(D_i)\subset {\rm int}\, D_i$, $i=1,2$ (see Fig. \ref{split}). 
\begin{figure}[h!]
			\centerline{\includegraphics [width=8 cm]{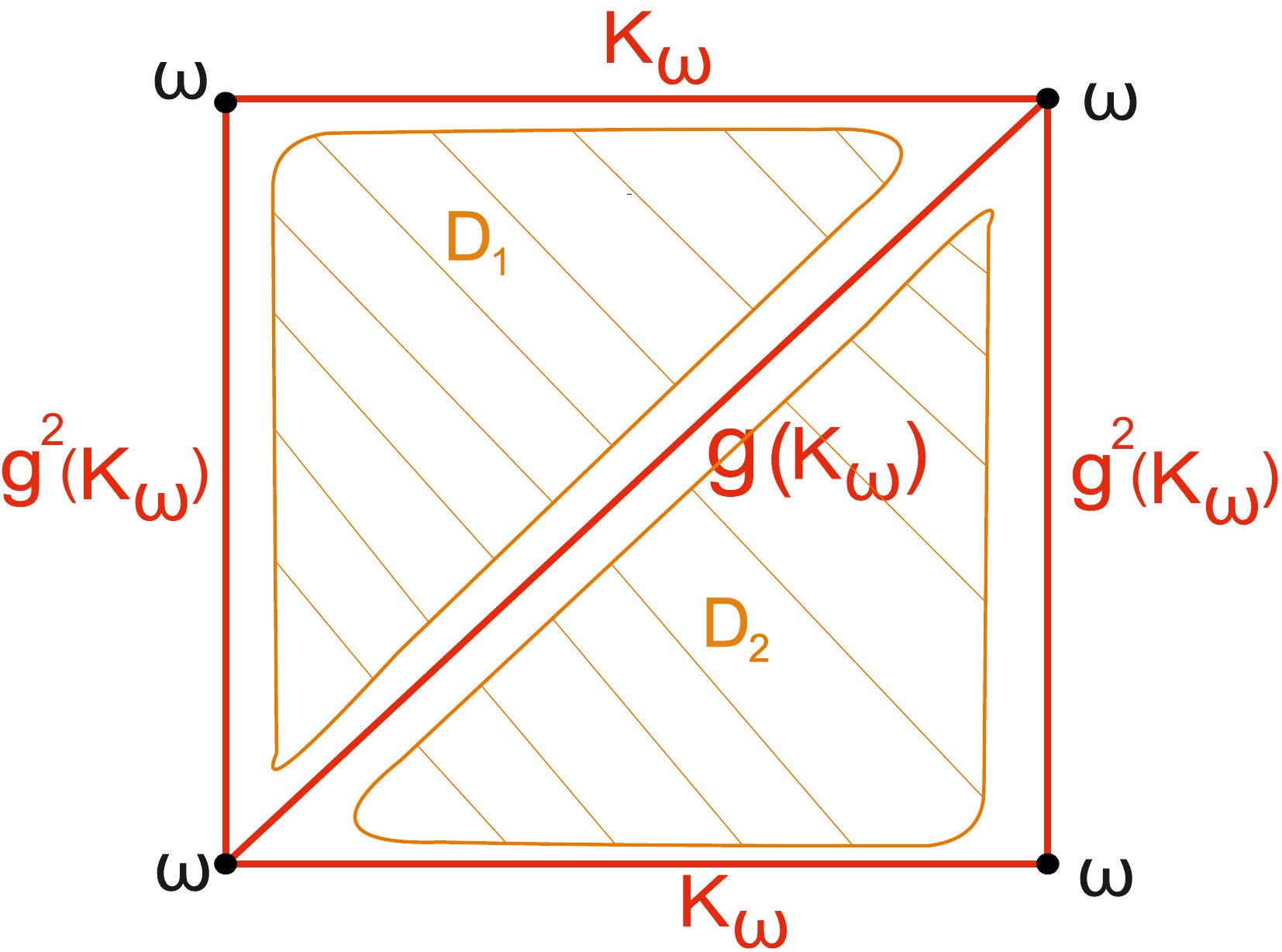}}
			\caption{\small Splitting of $\mathbb T^2$ by $A_\omega$}\label{split}
		\end{figure} 
By the Brouwer fixed-point theorem, the map $g$ has at least one fixed point in $D_i$. As $g$ has exactly two fixed point, excepted $\omega$, then $\alpha_i\in D_i$.  Due to Statement \ref{sad-noo}, there exists a stable arc $\tilde{g}_t:\mathbb T^2 \to\mathbb T^2$ such that $\tilde{g}_0=g$ and $\tilde g_1$ is a gradient-like diffeomorphism with the following properties:  $\Omega_{\tilde g_1}|_{D_i}=\alpha_i$ and  ${g_1}|_{U_{A_{\omega}}}=g|_{U_{A_{\omega}}}$. It means that $\tilde g_1$ is the simplest.

2) In this case we show that diffeomorphism $g$ is connected by a stable arc $\tilde g_t$ with a gradient-like diffeomorphism $\tilde g_1$ whose non-wandering set contains  fewer saddle points than $g$ does. Therefore, there is a stable arc, connecting $g$ with a diffeomorphism satisfying to the case 1), that finishes the proof.

For this aim we consider two subcases: 2a) $\Gamma_{\omega}$ contains a contractible loop, 2b) $\Gamma_{\omega}$ does not contain a loop. 

2a) Denote by $K_\omega$ a contractible loop in $\Gamma_{\omega}$ and by $d_\omega\subset\mathbb T^2$ -- 2-disc with the boundary $K_\omega$. Then $A_\omega=d_\omega\cup g(d_\omega)\cup g^2(d_\omega)$ is an attractor of the diffeomorphism $g$ (see Fig. \ref{Dw}).
\begin{figure}[h!]
			\centerline{\includegraphics [width=8 cm]{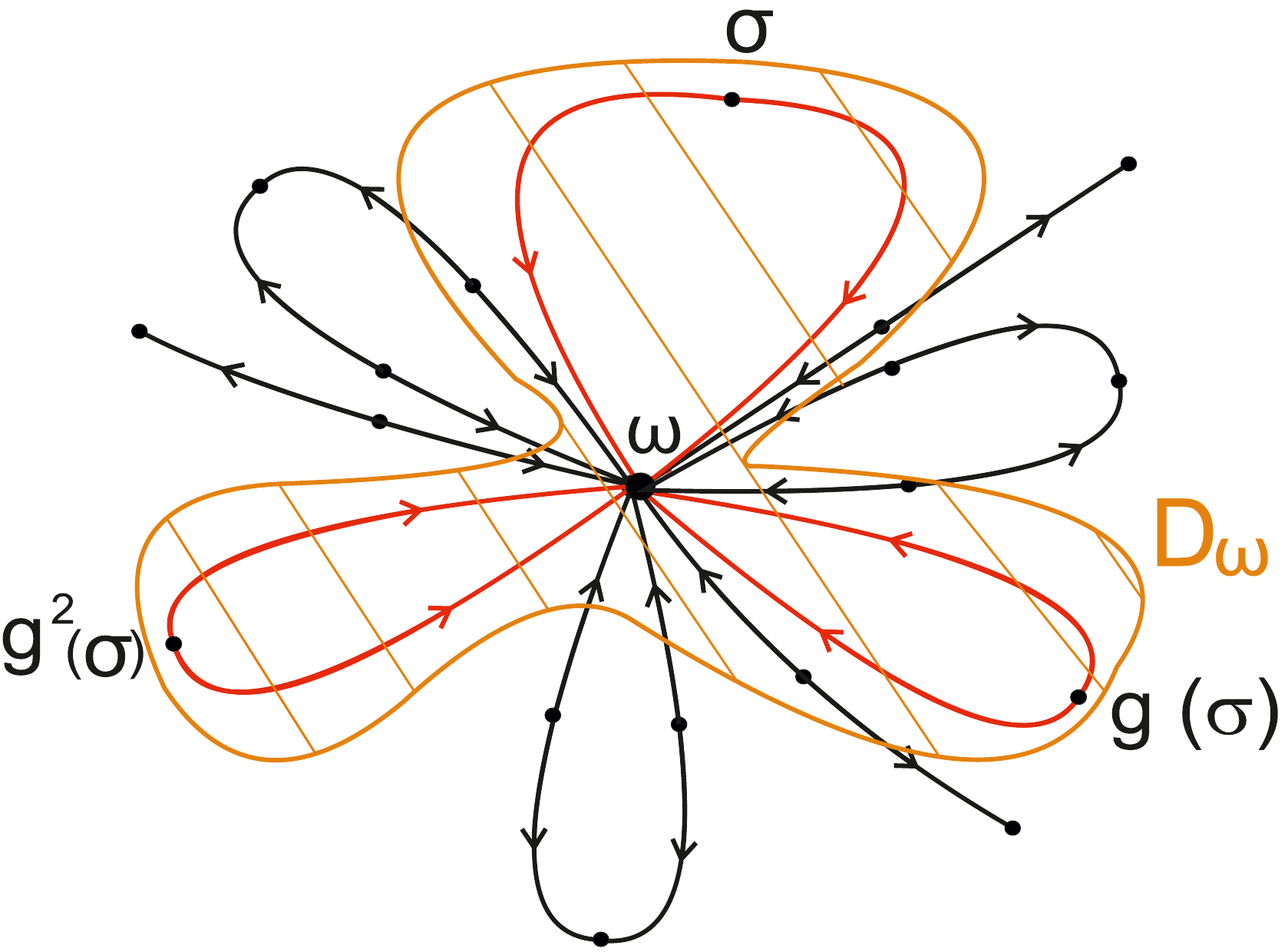}}
			\caption{\small Complex $\Gamma_\omega$ with contractible loop}\label{Dw}
		\end{figure}
Then there is its trapping  neighborhood, which is a 2-discs $D_\omega$. As $D_\omega$  contains a fixed sink $\omega$, then, by Statement \ref{sad-noo}, there exists a stable arc $\tilde{g}_t:\mathbb T^2 \to\mathbb T^2$ such that $\tilde{g}_0=g$ and $\tilde g_1$ is a gradient-like diffeomorphism with the following properties:  $\Omega_{\tilde g_1}|_{D_\omega}=\omega$ and  ${g_1}|_{\mathbb T^2\setminus D_{\omega}}=g|_{\mathbb T^2\setminus D_{\omega}}$. It means that $\tilde g_1$ contains  fewer saddle points than $g$ does.

2b) If $\Gamma_\omega$ has no cycles, then, by Statement \ref{<n-2}, is a an attractor of the diffeomorphism $g$ with a trapping 2-disc $D_\omega$ (see Fig. \ref{Sw}).
\begin{figure}[h!]
			\centerline{\includegraphics [width=6 cm]{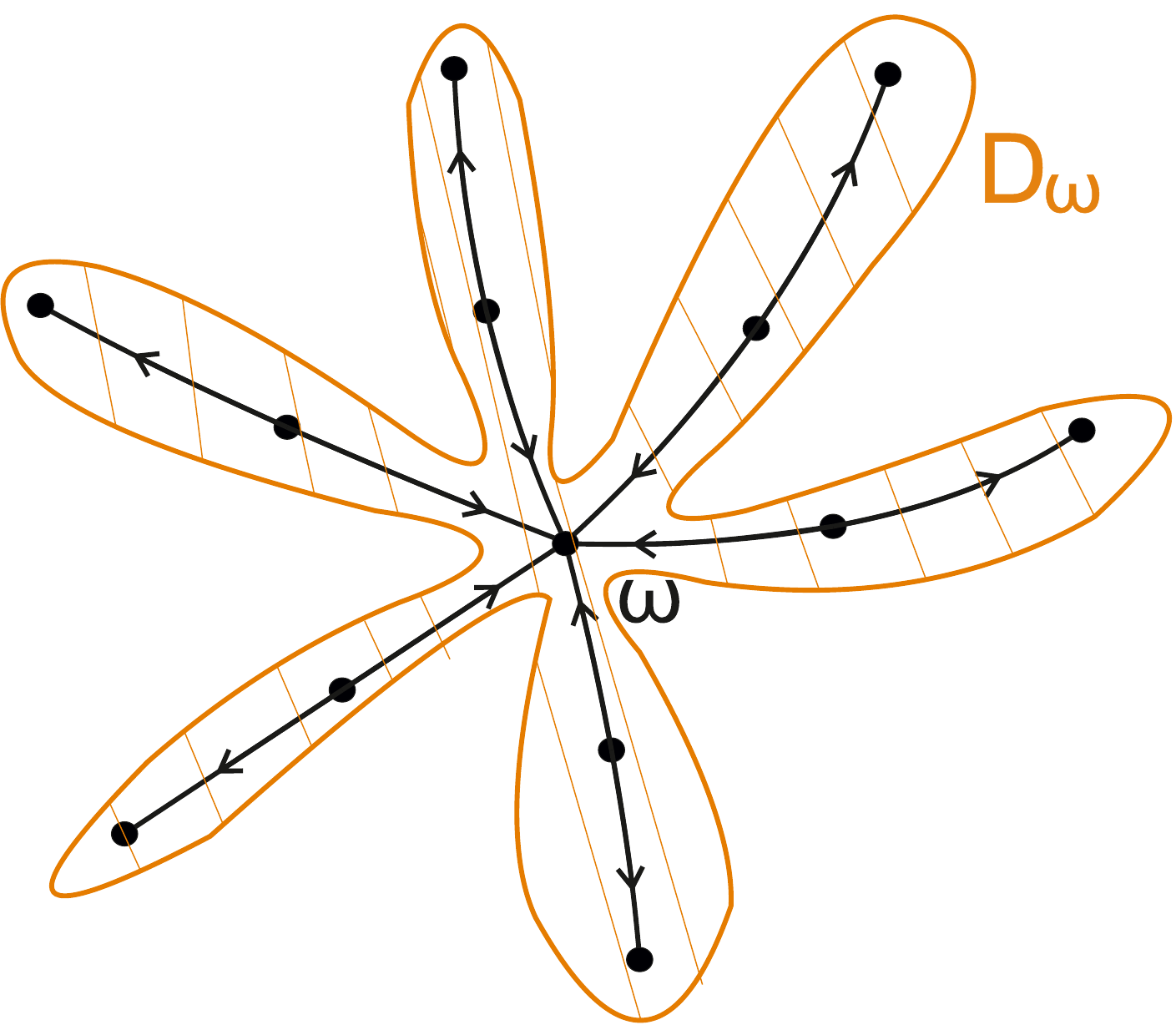}}
			\caption{\small Complex $\Gamma_\omega$ without loop}\label{Sw}
		\end{figure}
Similarly to the case 2a) it is possible to show that $g$ is connected by a stable arc $\tilde g_t$ with a gradient-like diffeomorphism $\tilde g_1$ whose non-wandering set contains  fewer saddle points than $g$ does.

{\bf $i=0$.} Let $g\in G_{2,0}$. Firstly, notice, that by Lemma \ref{dyng}, all saddle and sink points of the diffeomorphism $g$ have period 3. Then, by Statement \ref{sad-noo}, we can cancel all saddle points whose unstable manifolds contains in its closure sink points from different orbits (see Fig. \ref{canc}), because such a structure is in the $2$-disk.  
\begin{figure}[h!]
			\centerline{\includegraphics [width=12 cm]{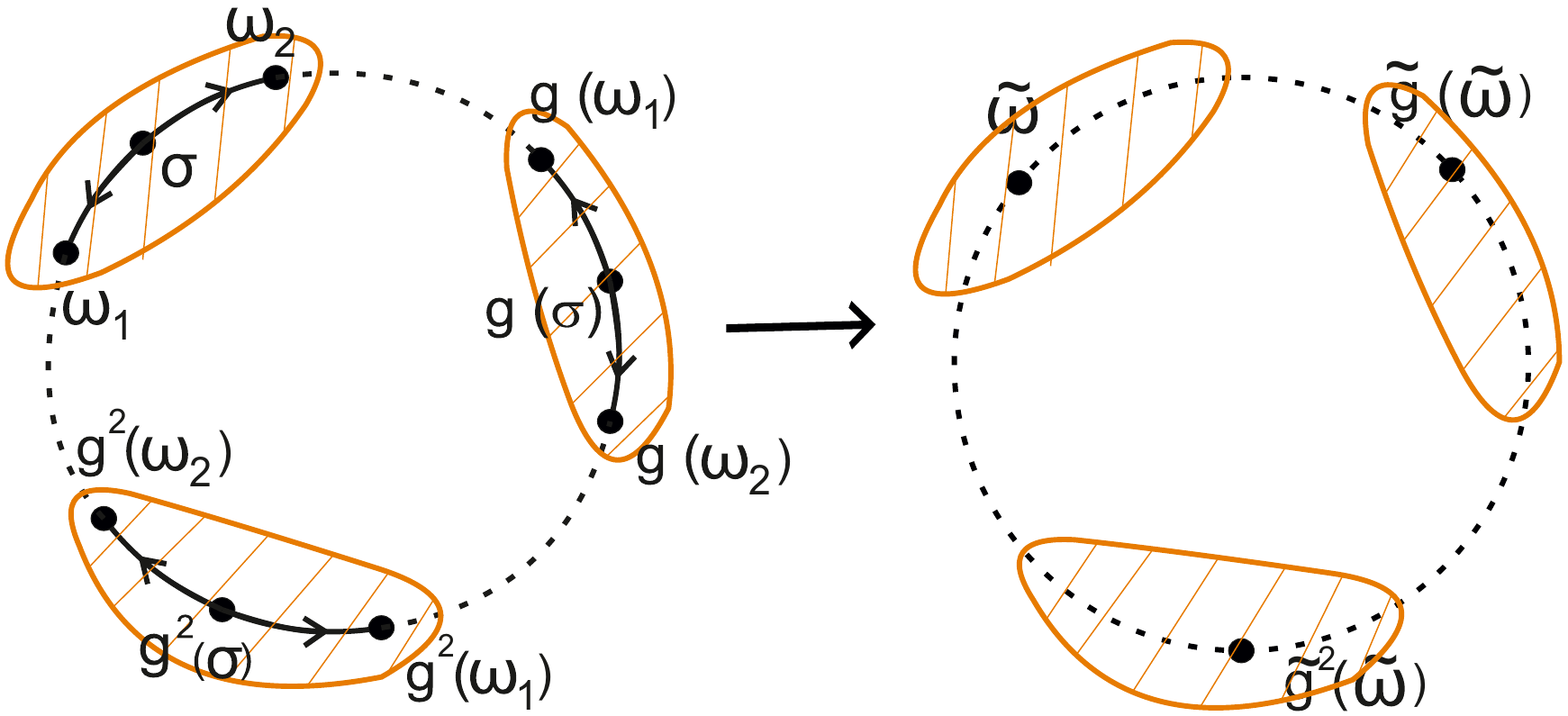}}
			\caption{\small Cancelling of saddle points with the different unstable asymptotic}\label{canc}
		\end{figure}
Also, we can reduce all saddle points $\sigma$ whose unstable manifolds contains in its closure a unique sink point $\omega$. Because knots ${\rm cl}\,W^u_\sigma, g({\rm cl}\,W^u_\sigma),g^2({\rm cl}\,W^u_\sigma)$ are pairwise disjoint and, hence, have the same homotopic type (see, for example, \cite{Ko}) which preserves by $g_\bigstar=A_2$. Such type is unique and it is $\langle 0,0\rangle$, that is the knots bound 2-discs. By Statement \ref{sad-noo} we can change it (along a stable arc) by a sink orbit (see Fig. \ref{loo}). 
\begin{figure}[h!]
			\centerline{\includegraphics [width=11 cm]{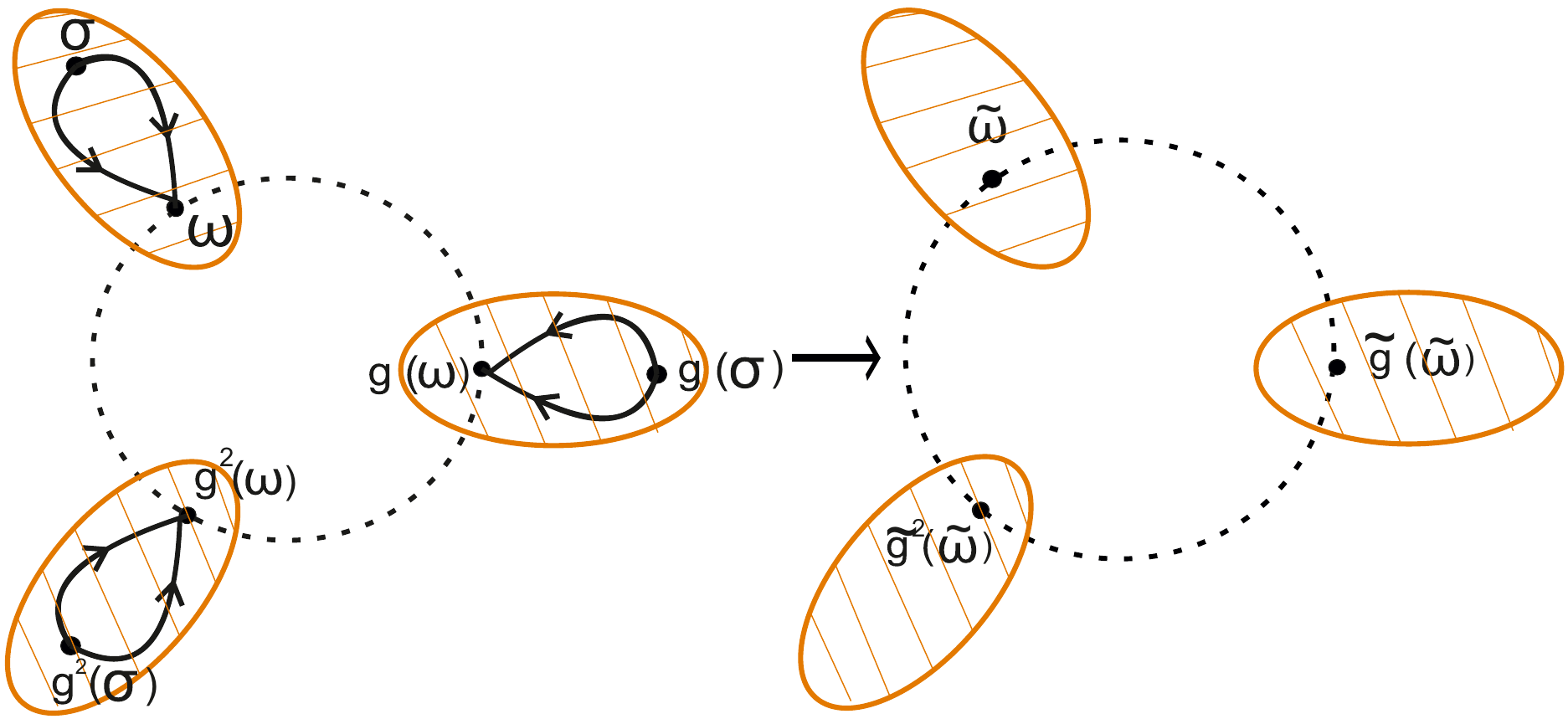}}
			\caption{\small Cancelling of saddle points with the same unstable asymptotic}\label{loo}
		\end{figure}

Thus, we can assume that any saddle point of the diffeomorphism $g$ have in its unstable manifold closure sinks from the same orbit. On the other side, the set ${\rm cl}\,W^u_{\Omega^1_g}$ contains all sink points and, by Statement \ref{<n-2}, it is connected.  
So, $\Omega^0_g$ consists of the unique sink orbit $\omega,g(\omega),g^2(\omega)$. By  Statement \ref{lefschetz} the set $\Omega^1_g$ contains at lest two points. It means that there are different saddle points $\sigma_1,\sigma_2$ such that $${\rm cl}\,W^u_{\sigma_1}\setminus W^u_{\sigma_1}={\rm cl}\,W^u_{\sigma_2}\setminus W^u_{\sigma_2}=\omega\sqcup g(\omega).$$  Then the set $K=\omega\cup W^u_{\sigma_1}\cup W^u_{\sigma_2}\cup g_0(\omega)$ is a knot on the torus. If the knot is contractible, we can confluence orbits of $\sigma_1$ and $\sigma_2$ (see Fig. \ref{conf}) by Statement \ref{ann}, and choose after that two new saddle points $\sigma_1,\sigma_2$.  So, we assume that $K$ is non-contractible.
\begin{figure}[h!]
			\centerline{\includegraphics [width=11 cm]{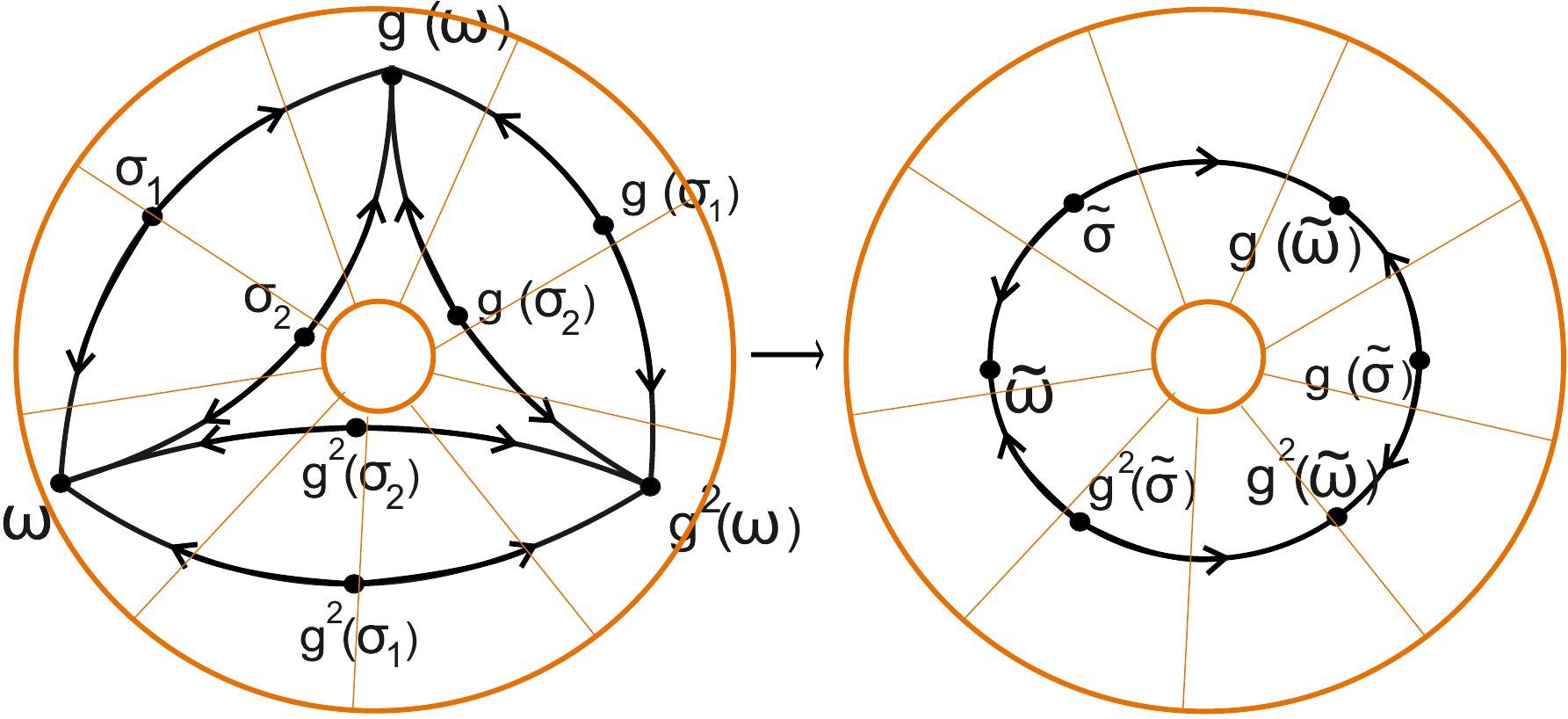}}
			\caption{\small Confluence of saddle orbits}\label{conf}
		\end{figure} 
		
Since $g_{\bigstar}=A_2$, similarly to item 3 of Lemma \ref{g1}, it is possible to prove that  the knots $K,g(K),g^2(K)$ have the following homotopy types: $\pm\langle 1,0\rangle$, $\mp\langle 1,1\rangle$, $\pm\langle 0,1\rangle$. Moreover, by Statement \ref{<n-2}, $$A=K\cup g(K)\cup g^2(K)$$ is an attractor of the diffeomorphism $g$. Then there is its trapping  neighborhood $U_{A_{\omega}}$ such that  $\mathbb T^2 \setminus{\rm int}\,U_{A_{\omega}}$ is the disjoint union of two 3-discs $D_1 \sqcup D_2\sqcup$ and $g^{-1}(D_i)\subset {\rm int}\, D_i$, $i=1,2,3$ (see Fig. \ref{split0}). 
\begin{figure}[h!]
			\centerline{\includegraphics [width=6 cm]{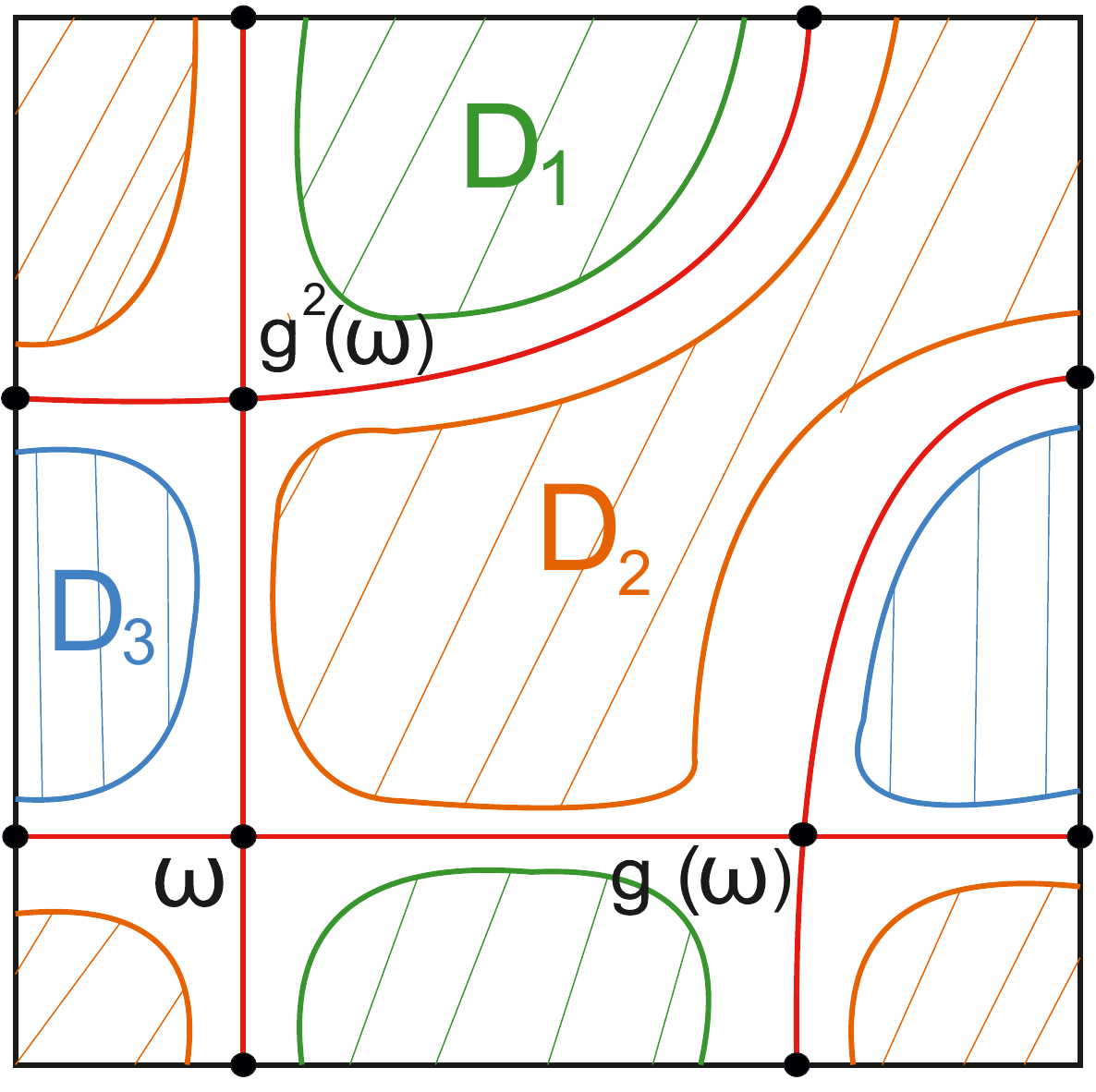}}
			\caption{\small Splitting of $\mathbb T^2$ by $A$}\label{split0}
		\end{figure} 
By the Brouwer fixed-point theorem, the map $g$ has at least one fixed point in $D_i$. As $g$ has exactly three fixed point, then $\alpha_i\in D_i$.  Due to Statement \ref{sad-noo}, there exists a stable arc $\tilde{g}_t:\mathbb T^2 \to\mathbb T^2$ such that $\tilde{g}_0=g$ and $\tilde g_1$ is a gradient-like diffeomorphism with the following properties:  $\Omega_{\tilde g_1}|_{D_i}=\alpha_i$ and  ${g_1}|_{U_{A_{\omega}}}=g|_{U_{A_{\omega}}}$. It means that $\tilde g_1$ is the simplest.
\end{proof}
	
	\begin{lemma}\label{step2} Any the simplest diffeomorphisms of class $G_{2,i}$ are connected by an arc without bifurcations.
	\end{lemma}
	\begin{proof} We will carry out the proof for the simplest diffeomorphisms $g_1,g'_1\in G_{2,1}$, for the remaining classes $G_{2,i}$ the proof is carried out similarly.
		
		According to Lemma \ref{g1}, the union of the closures of all unstable manifolds of saddle points of diffeomorphisms $g_1,g'_1$ consists of three knots $K,g_1(K),g_1^2(K);\,K',g'_1(K'),g_1^{\prime 2}(K')$, respectively. In this case, the knots of one family are pairwise homotopic to the knots from another family. Then (see, for example, \cite{Rol}) there exists a diffeotopy $h_t:\mathbb T^2\to\mathbb T^2$ such that $h_0=id$ and the closures of the unstable saddle manifolds of the diffeomorphism $\tilde g_1=h_1g_1h_1^{-1}$ coincide with similar manifolds of the diffeomorphism $g_1$. According to \cite[Lemmas 7.2,7.3]{NoPo22}, the diffeomorphism $\tilde g_1$ is connected by an arc without bifurcations to the diffeomorphism $g_1$.\end{proof}

\section*{Author information}
\subsection*{Authors and Affiliations}
Denis Baranov, HSE University, Russia, dabaranov@hse.ru.\\
Olga Pochinka, HSE University, Russia, olga-pochinka@yandex.ru.

\subsection*{Corresponding author}
Correspondence to Olga Pochinka.

\section*{Ethics declarations}
\subsection*{Conflict of Interest} 
No potential conflict of interest was reported by the authors.

\section*{Funding}	

The work was supported by the Russian Science Foundation, project no. 22-11-00027-$\Pi$.

	\end{document}